\documentclass[11pt,reqno]{amsart}

\usepackage[unicode=true,pdfusetitle,
 bookmarks=true,bookmarksnumbered=false,bookmarksopen=false,
 breaklinks=false,pdfborder={0 0 0},backref=false,colorlinks=true]
 {hyperref}
 
\usepackage{color}
\usepackage{amsfonts,amssymb}
\usepackage{bbm} 
\usepackage{marginnote}
\usepackage{pdfsync}
\usepackage{a4wide}

\usepackage{graphicx} 
\usepackage{pstricks}
\usepackage{pinlabel}

\usepackage{xfrac}

\usepackage{varioref}

\numberwithin{equation}{section}

\newtheorem{neu}{}[section]
\newtheorem{Cor}[neu]{Corollary}
\newtheorem*{Cor*}{Corollary}
\newtheorem{Thm}[neu]{Theorem}
\newtheorem*{Thm*}{Theorem}
\newtheorem{Lemma}[neu]{Lemma}
\newtheorem*{Observation*}{Observation}
\newtheorem{Prop}[neu]{Proposition}
\newtheorem*{Prop*}{Proposition}
\theoremstyle{definition}
\newtheorem*{Rmk*}{Remark}
\newtheorem{Rmk}[neu]{Remark}
\newtheorem{Ex}[neu]{Example}
\newtheorem*{Ex*}{Example}

\newtheorem*{Qu*}{Question}

\newtheorem{Def}[neu]{Definition}

\theoremstyle{remark}
\newtheorem*{acknowledgement*}{\bf Acknowledgement}

\newcommand{\N}{\mathbbm{N}}

\newcommand{\Z}{\mathbbm{Z}}
\newcommand{\R}{\mathbbm{R}}

\newcommand{\id}{\mathrm{id}}

\newcommand{\Cont}{\mathrm{Cont}}

\newcommand{\HF}{\mathrm{HF}}

\newcommand{\RFH}{\mathrm{RFH}}

\newcommand{\HM}{\mathrm{HM}}

\begin{document}
\title{Orderability, contact non-squeezing, and Rabinowitz Floer homology}
\author{Peter Albers}
\author{Will J.~Merry}
\address{
	Peter Albers\\ 
	Mathematisches Institut\\
	Westf\"alische Wilhelms-Universit\"at M\"unster}
\email{peter.albers@wwu.de}
\address{
 	Will J.~Merry\\
	Department of Mathematics\\
	ETH Z\"urich}
\email{merry@math.ethz.ch}

\begin{abstract}
We study Liouville fillable contact manifolds $(\Sigma,\xi)$ with non-zero Rabinowitz Floer homology and assign spectral numbers to paths of contactomorphisms. As a consequence we prove that $\widetilde{\mathrm{Cont}_0}(\Sigma,\xi)$ is orderable in the sense of Eliashberg and Polterovich. This provides a new class of orderable contact manifolds. If the contact manifold is in addition periodic or a prequantization space $M \times S^1$ for $M$ a Liouville manifold, then we construct a contact capacity. This can be used to prove a general non-squeezing result, which amongst other examples in particular recovers the beautiful non-squeezing results from \cite{EliashbergKimPolterovich2006}.
\end{abstract}
\maketitle

\section{Introduction and Results}
\label{sec:Introduction}

Suppose $(\Sigma,\xi)$ is a closed coorientable contact manifold.
Denote by $\mathrm{Cont}_{0}(\Sigma,\xi)$ the identity component of
the group of contactomorphisms, and denote by $\mathcal{P}\mathrm{Cont}_{0}(\Sigma,\xi)$
the set of smooth paths of contactomorphisms starting at the identity.
The universal cover $\widetilde{\mathrm{Cont}}_{0}(\Sigma,\xi)$ is
then $\mathcal{P}\mathrm{Cont}_{0}(\Sigma,\xi)/\sim$, where $\sim$
denotes the equivalence relation of being homotopic with fixed endpoints.
Suppose $\alpha\in\Omega^{1}(\Sigma)$ is a contact form defining
$\xi$, and $\theta^{t}$ its Reeb flow. To a path $\varphi=\{\varphi_{t}\}_{0\leq t\leq1}\in\mathcal{P}\mathrm{Cont}_{0}(\Sigma,\xi)$
we can associate its contact Hamiltonian $h_{t}$ 
\begin{equation}
h_{t}\circ\varphi_{t}:=\alpha\left(\frac{d}{dt}\varphi_{t}\right):\Sigma\rightarrow\mathbbm{R},
\end{equation}
which uniquely determines the path $\varphi$. In this article we
are interested in four classes of contact manifolds, labelled
\textbf{(A)}, \textbf{(A)}$^+$, \textbf{(B)}, and \textbf{(C)}. See Section \ref{sec:Preliminaries} for precise definitions of the terms
involved. 

\begin{description}
\item [Assumption {(A)}] $(\Sigma,\xi)$ admits a Liouville filling $W$ such that
the Rabinowitz Floer homology $\RFH_{*}(\Sigma,W)$ is non-zero. 
\end{description}

\begin{Thm}
\label{Theorem-A}Suppose $(\Sigma,\xi)$ satisfies Assumption \textbf{\emph{(A)}}.
Then for any non-zero class $Z\in\RFH_{*}(\Sigma,W)$ there is a map $c(\cdot,Z):\mathcal{P}\mathrm{Cont}_{0}(\Sigma,\xi)\rightarrow\mathbbm{R}$ with the following properties.
\begin{enumerate}\itemsep0.5ex
\item If $\varphi\sim\psi$ then $c(\varphi,Z)=c(\psi,Z)$. Thus $c(\cdot,Z)$ descends
to define a map (denoted by the same symbol) $c(\cdot,Z):\widetilde{\mathrm{Cont}}_{0}(\Sigma,\xi)\rightarrow\mathbbm{R}$.
\item For any $T\in\mathbbm{R}$, $c(t\mapsto\theta^{tT},Z)=-T+c(\id_\Sigma,Z)$.
\item The map $c$ is continuous with respect to the $C^{2}$-norm
on $\mathcal{P}\mathrm{Cont}_{0}(\Sigma,\xi)$.
\item If $\varphi$ resp.~$\psi$ is generated by the contact Hamiltonian $h_{t}$ resp.~$k_t$ with
$h_t(x)\geq k_t(x)$ for all $x\in\Sigma$ and $t\in[0,1]$ then $c(\varphi,Z)\leq c(\psi,Z)$. 
\end{enumerate}
\end{Thm}

\begin{Cor}\label{cor:positive_implies_negative}
If there exists a constant $\delta>0$ such that $h_t(x)\geq \delta$ for all $x\in\Sigma$ and $t\in[0,1]$ then $c(\varphi,Z)<c(\id_\Sigma,Z)$.
\end{Cor}

\begin{Cor}\label{cor:(A)_implies_orderable}
If $(\Sigma,\xi)$ satisfies Assumption \textbf{\emph{(A)}} then $\widetilde{\mathrm{Cont}}_{0}(\Sigma,\xi)$ is orderable in the sense of Eliashberg--Polterovich \cite{EliashbergPolterovich2000}.
\end{Cor}

\begin{proof}[Proof of Corollary \ref{cor:positive_implies_negative}]
Note that the constant function $\delta$ generates the path $\{t\mapsto\theta^{t\delta}\}$ thus Theorem \ref{Theorem-A} (2) \& (4) implies
\begin{equation}
c(\varphi,Z)\leq c(t\mapsto\theta^{t\delta},Z)=-\delta+c(\id_\Sigma,Z)<c(\id_\Sigma,Z).
 \end{equation}
\end{proof}

\begin{proof}[Proof of Corollary \ref{cor:(A)_implies_orderable}]
Recall from \cite[Criterion 1.2.C.]{EliashbergPolterovich2000} that $\widetilde{\mathrm{Cont}}_{0}(\Sigma,\xi)$ is orderable if and only if no contractible loop $\varphi$ of contactomorphisms exists whose contact Hamiltonian $h_t$ satisfies $h_t(x)>0$ for all $x\in\Sigma$ and $t\in[0,1]$. Let us assume, by contradiction, that $\varphi$ is such a loop. Then (1) in Theorem \ref{Theorem-A} implies that $c(\varphi,Z)=c(\id_\Sigma,Z)$ since $\varphi$ is contractible. On the other hand Corollary \ref{cor:positive_implies_negative} implies that $c(\varphi,Z)<c(\id_\Sigma,Z)$. This contradiction proves the Corollary.
\end{proof}

\begin{Rmk}
Together with its companion article \cite{AlbersFuchsMerry2013} this article is the first to establish Rabinowitz Floer homology as a tool for studying orderability and non-squeezing questions in contact geometry. The aim of the article \cite{AlbersFuchsMerry2013} is very different from the present one since it is solely concerned with a link between the famous Weinstein conjecture and orderability. In this article we derive obstructions from Rabinowitz Floer homology to non-orderability and to squeezing phenomena. Since Rabinowitz Floer homology is nowadays rather computable this delivers checkable criteria for orderability and non-squeezing. In particular, we reproduce many of the previously known examples of orderable contact manifolds and similarly for the non-squeezing results. At the same time our approach gives entirely new classes of orderable contact manifolds and an abstract non-squeezing results.

A precursor to this development is the article by Frauenfelder \cite{AlbersFrauenfelder2012} and the first author in which a rather different version of Rabinowitz Floer homology is used to mimic Givental's construction of the non-linear Maslov index. On unit cotangent bundles this also leads to an obstruction to a (strong form) of non-orderability. 
\end{Rmk}

Corollary \ref{cor:(A)_implies_orderable} has the following rephrasing.

\begin{Cor}
Let $(\Sigma,\xi)$ be a closed contact manifold for which $\widetilde{\mathrm{Cont}}_{0}(\Sigma,\xi)$ is not orderable. Then for any Liouville filling $W$ of $\Sigma$ one has
\begin{equation}
\RFH_{*}(\Sigma,W)=0.
\end{equation}
\end{Cor}

We illustrate the above at some examples.

\begin{Ex}$ $ 
\begin{itemize}
\item The sphere $S^{2n-1}$ with its \emph{standard} contact structure is not orderable by \cite[Theorem 1.10]{EliashbergKimPolterovich2006}. The equivalent statement of vanishing symplectic homology of any Liouville filling of the standard contact sphere was proved before by Smith, see \cite[Corollary 6.5]{Seidel2006}.
\item A new class of orderable contact manifolds is given by links of weighted homogeneous singularities with positive Milnor number. This includes certain \emph{Brieskorn manifolds}, and in particular \emph{non-standard} structures on spheres (the Ustilovsky spheres), as well as contact structures on \emph{exotic} spheres. This was communicated to us by Otto van Koert, see Example \ref{Ex:Otto} below, which also includes more examples.
\item  Let $\Sigma=S^*_gB$ be the unit cotangent bundle of the closed manifold $B$ equipped with its standard contact structure $\xi$. The Liouville filling by the unit codisk bundle $D^*_gB$ always has $\RFH_*(S^*_gB,D^*_gB)\neq0$ due to Cieliebak-Frauenfelder-Oancea \cite{CieliebakFrauenfelderOancea2010}, see also Abbondandolo-Schwarz \cite{AbbondandoloSchwarz2009}. Thus, $\widetilde{\mathrm{Cont}}_{0}(S^*_gB,\xi)$ is orderable, which was proved by Eliashberg-Kim-Polterovich \cite{EliashbergKimPolterovich2006} and Chernov-Nemirovski \cite{ChernovNemirovski2010}.
\item Since symplectic homology and thus Rabinowitz Floer homology is unchanged under subcritical handle attachment (\cite{Cieliebak2002}) any Liouville fillable contact manifold with non-zero Rabinowitz Floer homology remains $\widetilde{\mathrm{Cont}}_{0}$-orderable under subcritical handle attachment.
\end{itemize}
\end{Ex}

A slight strengthening of Assumption $\textbf{(A)}$ is the following, see the discussion \vpageref{ass:a+} for details.

\begin{description}
\item [Assumption {(A)$ ^+$}] $(\Sigma,\xi)$ admits a Liouville filling $(W,d\lambda)$
such that $\alpha:=\lambda|_{\Sigma}$ is Morse-Bott. Moreover there exists a non-zero class  $\mu_\Sigma \in \RFH_n( \Sigma,W)$ such that $  p^{ - \varepsilon}(  \mu_{ \Sigma} ) = j^{ \varepsilon}( [ \Sigma])$.
\end{description}

\begin{Ex}
\label{Ex:hypertight}
If  $(\Sigma,\xi)$ admits a Liouville filling $(W,d\lambda)$ such that the Reeb vector field of $\alpha:=\lambda|_{\Sigma}$ has no contractible Reeb orbits (e.g. $\mathbbm{T}^3$ with its standard contact structure, which is filled by $D^*\mathbbm{T}^2$) then Assumption $\textbf{(A)}^+$ is trivially satisfied, since in this case $\RFH_*(\Sigma ,W) \cong \mathrm{H}_{* + n - 1}(\Sigma ;\Z_2)$. 
\end{Ex}

The class $\mu_{\Sigma}$ has the property that $c( \id_{\Sigma}, \mu_{\Sigma}) = 0$.  Assumption $\textbf{(A)}^+$  allows us to strengthen Statement (4) of Theorem \ref{Theorem-A} to the following statement:

\begin{Thm}
\label{thm:weak_inequality}
Suppose $\varphi\in\widetilde{\mathrm{Cont}}_{0}(\Sigma,\xi)$
has contact Hamiltonian $h_{t}$. Assume $h_t\leq0$ and there exists $x\in\Sigma$ such that $h_t(x)<0$ for all $t\in[0,1]$ then $c(\varphi,\mu_\Sigma)>0$.
\end{Thm}

\begin{Rmk}
In Section \ref{sub:Prequantisation-spaces} below we provide an example to show that the same implication with opposite inequalities in the above theorem does not hold. See Remarks \ref{rem:reversing_prev_cor} and  \ref{rem:cptly_sprtd_rf}, and Appendix \ref{app:A-'compactly-supported}. 

It's worthwhile pointing out that replacing $\alpha$ with $-\alpha$ changes the
sign of the contact Hamiltonian $h_{t}$. However in general if $(\Sigma,\alpha)$
is Liouville fillable then there is no canonical Liouville filling
of $(\Sigma,-\alpha)$, and thus for our purposes the distinction
between positive and negative is not arbitrary. A good non-trivial
example to bear in mind is $\Sigma=\mathbbm{R}P^{3}$, with $\xi$
the contact structure induced from from the standard contact structure
on $S^{3}$. Then $\Sigma$ is Liouville fillable with $W=T^{*}S^{2}$,
and $\RFH_{*}(\mathbbm{R}P^{3},T^{*}S^{2})$ is infinite dimensional
(cf. \cite{AbbondandoloSchwarz2009,CieliebakFrauenfelderOancea2010}).
However we can also view $\mathbbm{R}P^{3}$ as a prequantization space
over $S^{2}$, which corresponds to replacing $\alpha$ with $-\alpha$.
Indeed, in this case the completion $W$ is the total space of $\mathcal{O}(2)\rightarrow S^{2}$,
as opposed to $T^{*}S^{2}$ which is the total space of $\mathcal{O}(-2)\rightarrow S^{2}$. 

For us the main relevance of Theorem \ref{thm:weak_inequality} is that it implies the contact capacity $\bar{c}( \cdot, \mu_{\Sigma})$ we define below is  non-trivial, see Remark \ref{rmk:capacity_non_trivial} below.
\end{Rmk}

\begin{Def}
\label{def:periodic}
We call a contact form $\alpha$ periodic if its Reeb flow $\theta^t$ is a 1-periodic loop: $\theta^1=\id_\Sigma$.
\end{Def}

Let us now assume that $(\Sigma,\xi)$ satisfies the following condition:

\begin{description}
\item [Assumption {(B)}] $(\Sigma,\xi)$ admits a Liouville filling $(W,d\lambda)$
such that the Rabinowitz Floer homology $\RFH_{*}(\Sigma,W)$
is non-zero and such that $\alpha:=\lambda|_{\Sigma}$ is periodic. 
\end{description}

\begin{Ex}[Communicated to us by Otto van Koert]
\label{Ex:Otto}
An interesting class of examples where our results apply is the following. Assume that $(Q, \omega)$ is a simply connected symplectic manifold. Assume in addition that $[\omega]$ is integral and $(Q,\omega)$ is monotone, with non-positive monotonicity constant. Let $K \subset Q$ denote a closed connected symplectic submanifold of codimension 2 such that $K$ is Poincar\'e dual to $k [\omega]$ for some $k \in \N$. Such a hypersurface is known as a \emph{Donaldson hypersurface}, since Donaldson showed that every symplectic manifold with an integral symplectic form admits a symplectic submanifold Poincar\'e dual to $k [\omega]$ for $k \in \N$ sufficently large \cite{Donaldson1996}. Assume in addition that $\mathrm{H}_1(K ;\Z) = 0$. Let $\nu(K)$ denote a collar neighborhood of $K$ in $Q$. Then the complement $Q \setminus \nu(K)$ is the interior of a Liouville domain $(W_1 , \lambda_1)$ with the property that the Reeb flow on $\Sigma : = \partial W_1$ is periodic (see for instance \cite{Diogo2012}). Denote by $W$ the completion of $W_1$. If we assume that the inclusion $\Sigma \hookrightarrow W$ induces an injection on $\pi_1$ (e.g.~if $ \dim Q \ge 6$), then the symplectic homology $\mathrm{SH}_*(W)$ is non-zero (see below), and hence so is the Rabinowitz Floer homology $\RFH_*( \Sigma ,W)$. Thus Assumption $\textbf{(B)}$ is satisfied. 

There are several ways to see that the symplectic homology of $\mathrm{SH}_*(W)$ is non-zero. The simplest one is an index argument, and goes as follows. Since $(Q,\omega)$ is monotone it follows easily that $c_1(TW)$ is torsion, and hence the Conley-Zehnder index is a well defined integer for contractible orbits. Next the monotonicity assumption and a suitable choice of Hamiltonian functions imply that the index of all contractible Reeb orbit is at most $n=\frac12\dim W$. However there is a well defined map $\mathrm{H}_{* + n}(W_1, \Sigma) \to \mathrm{SH}_*(W)$, and the image of the fundamental class has degree $n$. This class therefore remains non-zero in $\mathrm{SH}_*(W)$ due to index reasons. Alternatively, one can argue using $S^1$-equivariant symplectic homology: the proof of Lemma 7.6 in \cite{ChiangDingvanKoert2014} implies that $\mathrm{SH}^{S^1,+}_*(W)$ has no generators with large positive degree, since the index growth of non constant one periodic orbits is proportional to the (non-positive) monotonicity constant. However if $\mathrm{SH}_*(W)=0$ then work of Bourgeois-Oancea \cite{BourgeoisOancea2012} implies that $\mathrm{SH}_*^{S^1}(W)=0$. The Viterbo long exact sequence (see \cite[Lemma 4.8]{BourgeoisOancea2013a}) then implies $\mathrm{SH}^{S^1,+}(W) \cong \mathrm{H}_{*+n}(W, \Sigma) \otimes \mathrm{H}_*( \mathbbm{C} \mathrm{P}^{\infty} ; \Z)$, which has generators with arbitrary positive degree, which is a contradiction.

If we assume that $\pi_2(Q) = 0$ then it follows from the homotopy exact sequence of the fibration that all the Reeb orbits on $\Sigma$ are non-contractible. If in addition $\mathrm{H}_1(Q ;\Z) = 0$ and $\mathrm{H}_1(K ;\Z) = 0$ then the construction described above yields examples satisfying Assumption $\textbf{(A)}^+$ and $\textbf{(B)}$ (cf. Example \ref{Ex:hypertight}). Moreover, very explicit examples are the complement of a degree $k$-curve in $\mathbbm{C} \mathrm{P}^2$ with $k\geq 3$ which admit Liouville fillings (even though $\mathrm{H}_1(K ;\Z) \neq 0$).

Finally, another more general class of examples where our results apply are links of weighted homogeneous singularities with positive Milnor number. The latter guarantees the existence of Lagrangian spheres which in turn implies non-vanishing of $\RFH$, see \cite{KwonvanKoert2013}. In particular, this includes certain \emph{Brieskorn manifolds}, see \cite[Theorem 1.2]{KwonvanKoert2013} for a precise statement. We mention here only that these include \emph{non-standard} contact structures on spheres (the Ustilovsky spheres \cite{Ustilovsky1999}), as well as contact structures on \emph{exotic} spheres.
\end{Ex}

The advantage of Assumption $\textbf{(B)}$ is the following. As before $Z$ denotes a non-zero class in $\RFH_{*}(\Sigma,W)$. The basic definitions and results that follow are based on Sandon's article \cite{Sandon2011b}.

\begin{Def}
\label{Def of c bar when resonant}We define for $\varphi\in\widetilde{\Cont}_{0}(\Sigma,\xi)$
an integer $\overline{c}(\varphi,Z)$ by 
\begin{equation}
\overline{c}(\varphi,Z):=\left\lceil c(\varphi,Z)\right\rceil .
\end{equation}
 \end{Def}

\begin{Prop}
\label{prop:c_conjugation_invariant_introduction}
The function $\overline{c}(\cdot,Z):\widetilde{\mathrm{Cont}}_{0}(\Sigma,\xi)\rightarrow\mathbbm{Z}$
is conjugation invariant: if $\psi\in\mathrm{Cont}_{0}(\Sigma,\xi)$
and $\varphi\in\widetilde{\mathrm{Cont}}_{0}(\Sigma,\xi)$ then
\begin{equation}
\bar{c}(\psi\varphi\psi^{-1},Z)=\bar{c}(\varphi,Z).
\end{equation}
\end{Prop}
\begin{Rmk}
\label{Rmk:contact_conjguation_invariance_not_obvious}
In contrast to spectral invariants in Hamiltonian Floer homology, Proposition \ref{prop:c_conjugation_invariant_introduction} is a non-trivial result. See Remark \ref{Rmk:hamiltonian_trivial} below.
\end{Rmk}
\begin{Def}
For an open set $U\subset\Sigma$ we define the \emph{contact capacity}
\begin{equation}
\overline{c}(U,Z):=\sup\left\{ \overline{c}(\varphi,Z)\mid\varphi\in\widetilde{\Cont}_{0}(\Sigma,\xi),\ \mathfrak{S}(\varphi)\subset U\right\} \in\mathbbm{Z}\cup\{ \pm \infty\},
\end{equation}
where by convention we declare that $ \sup \emptyset  = - \infty$.
\end{Def}

\begin{Rmk}
The notion of contact capacity was introduced by Sandon in \cite{Sandon2011b}. She was the first to discover a connection between translated points and orderability and other contact rigidity phenomena. 
\end{Rmk}

Proposition \ref{prop:c_conjugation_invariant_introduction} has the following immediate corollary.

\begin{Cor}
\label{prop:cu_co_are_invariant_under_contactos_and_monotone_wrt_inclusions_INTRODUCTION}
For all $\psi\in\mathrm{Cont}_{0}(\Sigma,\xi)$, one has
\begin{equation}
\overline{c}(\psi(U),Z)=\overline{c}(U,Z).
\end{equation}
\end{Cor}

Using the contact capacity we obtain the following abstract non-squeezing results.
 
\begin{Thm}
\label{cor:non_squeezing_INTRODUCTION}
Let $U\subset V\subset\Sigma$ be open sets and
assume that there exists $\varphi\in\mathrm{Cont}_{0}(\Sigma,\xi)$
with $\varphi(V)\subset U$. Then 
\begin{equation}
\overline{c}(U,Z)=\overline{c}(V,Z).
\end{equation}
In particular, if $\overline{c}(U,Z)<\overline{c}(V,Z)$ then there exists
no contact isotopy mapping $V$ into $U$. 
\end{Thm}

\begin{proof}
Suppose $\varphi$ is as in the statement of the Theorem. Then trivially we have
$\overline{c}(U,Z)\leq\overline{c}(V,Z)$ and $\overline{c}(\varphi(V),Z)\leq\overline{c}(U,Z)$.
By Proposition \ref{prop:cu_co_are_invariant_under_contactos_and_monotone_wrt_inclusions} we also have $\overline{c}(\varphi(V),Z)=\overline{c}(V,Z)$,
and hence $\overline{c}(U,Z)=\overline{c}(V,Z)$ as claimed.
\end{proof}

\begin{Rmk}
\label{rmk:capacity_non_trivial}
If we assume both Assumptions  $\textbf{(A)}^+$ and $\textbf{(B)}$ then it follows from Theorem \ref{thm:weak_inequality} that $\bar{c}(U, \mu_{\Sigma})  \ge 1$ for every non-empty open subset $U \subset \Sigma$, and hence the capacity $\bar{c}(\cdot, \mu_{\Sigma})$ is non-trivial. In the next section we provide a class of examples where $\bar{c}(\cdot, \mu_{\Sigma})$ is computable, and thus derive applications of Theorem \ref{cor:non_squeezing_INTRODUCTION}.
\end{Rmk}

\subsection{\label{sub:Prequantisation-spaces}Prequantization spaces }
Fix a Liouville manifold $(M,d\gamma)$\textit{\emph{ (i.e.~the completion
of Liouville domain, cf. Definition \ref{A-Liouville-manifold}).
The }}\textit{prequantization space }\textit{\emph{of $M$ is the
contact manifold $\Sigma:=M\times S^{1}$, equipped with the contact
structure $\xi:=\ker\,\alpha$, where}}
\begin{equation}
\alpha:=\gamma+d\tau,
\end{equation}
where $\tau$ is the coordinate on $S^{1}=\mathbbm{R}/\mathbbm{Z}$.
These contact manifolds are the last type we study in this paper:
\begin{description}
\item [Assumption {(C)}] $(\Sigma,\xi=\ker\,\alpha)$ is a prequantization space $\Sigma=M\times S^{1}$,
where $(M,d\gamma)$ is a Liouville manifold, and $\alpha=\gamma+d\tau$. 
\end{description}
Let $P_1$ denote a torus with a small disc removed, so that $\partial P_1=S^{1}$.
Equip $P_1$ with an exact symplectic form $d\beta_1$ such that
$\beta_1|_{\partial P_1}=d\tau$. Let $(P,d\beta)$ denote the
completion of the Liouville domain $(P_1,d\beta_1)$, and consider
\[
W:=M\times P,
\]
equipped with the symplectic form $d\lambda$ where $\lambda:=\gamma+\beta$.
Even though $\Sigma$ is periodic, $W$ is \emph{not} a Liouville
filling of $\Sigma$, and in fact $\Sigma$ does not satisfy either
Assumptions \textbf{(A)} or \textbf{(B)} - for instance, $\Sigma$
is non-compact. Nevertheless, it is still possible to define the Rabinowitz
Floer homology $\RFH_{*}(\Sigma,W)$, and we prove that 
\begin{equation}
\label{eq:what rfh is}
\RFH_{*}(\Sigma,W)\cong\mathrm{HF}_{*}(M)\otimes\mbox{H}_{*}(S^{1};\mathbbm{Z}_{2}).
\end{equation}
Here
\begin{equation}
\label{eq:what hf is}
\mathrm{HF}_{*}(M)\cong\mbox{H}^{n-*}(M;\mathbbm{Z}_{2})
\end{equation}
 denotes the Hamiltonian Floer homology of $M$, defined using \emph{compactly
supported }Hamiltonians (see Frauenfelder-Schlenk \cite{FrauenfelderSchlenk2007}). Moreover the Rabinowitz Floer homology $\RFH_*( \Sigma,W)$ constructed in this way satisfies the analogue of Assumption \textbf{(A)$^+$}, that is, there is a suitable non-zero class $ \mu_{ \Sigma} \in \RFH_n( \Sigma, W)$. Indeed, in this case one simply takes $ \mu_{ \Sigma}$ to be the image of the class $ \{ \mathrm{pt} \} \times [S^1] \in \mathrm{H}^0 (M ; \Z_2) \otimes\mathrm{H}_1(S^1; \Z_2)$ under the isomorphisms \eqref{eq:what rfh is} and \eqref{eq:what hf is}.

Since the Hamiltonian Floer homology is non-zero, one can associate
a spectral number $c_{M}(f)$ to any $f\in\widetilde{\mbox{Ham}}_{c}(M,d\gamma)$,
the universal cover of the group of compactly supported Hamiltonian
diffeomorphisms (see eg. Schwarz \cite{Schwarz2000} or Frauenfelder-Schlenk
\cite{FrauenfelderSchlenk2007}) . As in the contact case described
above, $c_{M}$ can then be used to define a \emph{symplectic capacity
}$c_{M}(\mathcal{O})$ for $\mathcal{O}\subset M$ open, by setting
\begin{equation}
c_{M}(\mathcal{O}):=\sup\left\{ c_{M}(f)\mid\mathfrak{S}(f)\subset\mathcal{O}\right\} .
\end{equation}

\begin{Rmk}
\label{Rmk:hamiltonian_trivial}
In contrast to the contact case (see Proposition \ref{prop:c_conjugation_invariant_introduction} and Remark \ref{Rmk:contact_conjguation_invariance_not_obvious}), the proof that $c_{M}(f)$ is invariant under conjugation, that is,
$c_{M}(hfh^{-1})=c_{M}(f)$ for $f\in\widetilde{\mbox{Ham}}_{c}(M,d\gamma)$
and $h\in\mbox{Symp}_{c}(M,d\gamma)$ is immediate, since in this
case the action spectrum of $hfh^{-1}$ is the same as the action
spectrum of $f$ (see for instance \cite[Chapter 5, Proposition 7]{HoferZehnder1994}).
This in turn immediately implies that $c_{M}$ is a symplectic capacity,
that is, $c_{M}(f(\mathcal{O}))=c_{M}(\mathcal{O})$ for any symplectomorphism
$f$ and any open set $\mathcal{O}\subset M$. 
\end{Rmk}
Going back to $\Sigma=M\times S^{1}$, let us denote by $\mathrm{Cont}_{0,c}(\Sigma,\xi)$
those contactomorphisms $\varphi$ with compact support. There is
a natural way to lift an element $f\in\widetilde{\mbox{Ham}}_{c}(M,d\gamma)$
to obtain an element $\varphi\in\widetilde{\mathrm{Cont}}_{0,c}(\Sigma,\xi)$,
as we now explain. The equation
\begin{equation}
f_{t}^{*}\gamma-\gamma=da_{t},\ \ \ a_{0}\equiv0,
\end{equation}
determines a smooth compactly supported function $a_{t}:M\to\mathbbm{R}$.
Define $\varphi_{t}:\Sigma\rightarrow\Sigma$ by
\begin{equation}
\varphi_{t}(y,\tau):=\bigl(f_{t}(y),\underset{\textrm{mod 1}}{\underbrace{\tau-a_{t}(y)}}\bigr).\label{eq:phi from f}
\end{equation}

As explained in Section \ref{sub:Hamiltonian-Floer-homology} below, one can define for any non-zero class $Z$ the spectral numbers $c(\varphi , Z )$ for $\varphi\in\widetilde{\mathrm{Cont}}_{0,c}(\Sigma,\xi)$ in much the same way as before. Similarly one can define the capacity $\bar{c}(U)$ for $U \subset \Sigma$ open in the same way as before (if $U$ is not precompact then one must again only use elements of  $\widetilde{\mathrm{Cont}}_{0,c}(\Sigma,\xi)$ when defining $\bar{c}(U)$). Moreover most of the results stated thus far in the paper continue to hold (this statement is made more precise in Section  \ref{sub:Hamiltonian-Floer-homology}). In particular, 
Parts (1), (3), and (4) of Theorem \ref{Theorem-A} remain true, and so do Proposition \ref{prop:c_conjugation_invariant_introduction} and Theorem
\ref{cor:non_squeezing_INTRODUCTION}.

 It is natural to ask the question: if $\varphi\in\widetilde{\mathrm{Cont}}_{0,c}(\Sigma,\xi)$
is the lift of $f\in\widetilde{\mbox{Ham}}_{c}(M,d\gamma)$, how is
$c(\varphi) := c( \varphi , \mu_ \Sigma) $ related to $c_{M}(f)$? Note that $\mathfrak{S}(\varphi)=\mathfrak{S}(f)\times S^{1}$,
and hence another question is how the $c_{M}$ capacity of $\mathcal{O}\subset M$
is related to the $c$ capacity of $\mathcal{O}\times S^{1}$. The
following result answers these questions in the nicest possible way. 
\begin{Thm}
\label{thm:equality of capacities}Suppose $f\in\widetilde{\mbox{\emph{Ham}}}_{c}(M,d\gamma)$,
and let $\varphi\in\widetilde{\mathrm{Cont}}_{0,c}(\Sigma,\xi)$ denote
the lift of $f$. Then 
\begin{equation}
c_{M}(f)=c(\varphi).
\end{equation}
Moreover, if $\mathcal{O}\subset M$ is open and has compact closure
then 
\begin{equation}
c_{M}(\mathcal{O})=c(\mathcal{O}\times S^{1}).
\end{equation}

\end{Thm}
Theorem \ref{thm:equality of capacities} allows us to prove non-squeezing
results on $\Sigma$ by making use of the known results on $M$. Let
$l^{t}:M\rightarrow M$ denote the flow of the Liouville vector field
on $M$ and set $\zeta^{t}:=l^{\log\, t}$. We will prove the following
general result.
\begin{Thm}
\label{thm:liouville flow}Suppose $\mathcal{O}\subset M$ is a non-empty
open set with compact closure and finite capacity: $c_{M}(\mathcal{O})<\infty$.
Suppose there exists a contact isotopy $\varphi\in\widetilde{\mathrm{Cont}}_{0,c}(\Sigma,\xi)$
such that 
\begin{equation}
\varphi_{1}\left(\zeta^{r_{2}}(\mathcal{O})\times S^{1}\right)\subset\zeta^{r_{1}}(\mathcal{O})\times S^{1}
\end{equation}
for $r_{1},r_{2}\in\mathbbm{R}$. Then $\left\lceil r_{2}\right\rceil \leq\left\lceil r_{1}\right\rceil $.
More generally if $\mathcal{O}\subset\mathcal{Q}\subset M$ are open
sets with the property that there exists $\varphi\in\widetilde{\mathrm{Cont}}_{0}(\Sigma,\xi)$
with $\varphi_{1}(\mathcal{Q}\times S^{1})\subset\mathcal{O}\times S^{1}$
then $\left\lceil c_{M}(\mathcal{Q})\right\rceil = \left\lceil c_{M}(\mathcal{O})\right\rceil $. 
\end{Thm}

\begin{proof}
Note that for any $r\in\mathbbm{R}$, 
\begin{equation}
c_{M}(\zeta^{r}(\mathcal{O}))=rc_{M}(\mathcal{O})\ne0,
\end{equation}
since $c_{M}(\mathcal{O})>0$ as $\mathcal{O}$ is non-empty. Thus
without loss of generality we may assume that $c_{M}(\mathcal{O})=1$.
Then 
\begin{equation}
\overline{c}(\zeta^{r}(\mathcal{O})\times S^{1})=\left\lceil c_{M}(\zeta^{r}(\mathcal{O}))\right\rceil =\left\lceil r\right\rceil .
\end{equation}
The result is now an immediate consequence of Theorem \ref{cor:non_squeezing_INTRODUCTION} (which, as remarked above,  does indeed remain true in this setting). The last statement follows similarly. 
\end{proof}
\begin{Rmk}
Theorem \ref{thm:liouville flow} recovers the beautiful non-squeezing
result of \cite[Theorem 1.2]{EliashbergKimPolterovich2006}. In this
case one takes $M=\mathbbm{R}^{2n}$ and $U$ the unit ball. They prove
that if $\left\lceil r_{1}\right\rceil <\left\lceil r_{2}\right\rceil $
then it is not possible to squeeze $B(r_{2})\times S^{1}$ into the
\emph{cylinder}%
\footnote{Note that whilst $C(r_{1}):=B^{2}(r_{1})\times\mathbbm{R}^{2n-2}$
does not have compact closure in $\mathbbm{R}^{2n}$, and thus $\overline{c}(C(r_{1})\times S^{1})$
is not defined, since we only work with compactly supported contactomorphisms
we can deduce this from the second statement of Theorem \ref{thm:liouville flow}
by taking $\mathcal{O}=B(r_{2})$ and $\mathcal{Q}$ a sufficiently
large ellipse contained in $C(r_{1})$. %
} $C(r_{1})\times S^{1}$. This result was also recovered by Sandon
\cite{Sandon2011b} using generating functions.
\end{Rmk}

A further applications of Theorem \ref{thm:equality of capacities}
is the following. Here we denote by $c_{\mathrm{HZ}}$ the \emph{Hofer-Zehnder capacity} (see Definition \ref{def_of_HZ} below or \cite{HoferZehnder1994}).
\begin{Thm}
\label{prop:our HZ result_INTRO
}Let $(M,d\gamma)$ denote a Liouville manifold.
Equip $\mathbbm{R}^{2m}$ with the standard symplectic form $d\lambda_{\mathrm{std}}$,
and consider the contact manifold $(\widetilde{\Sigma},\alpha+\lambda_{\mathrm{std}})$,
where $\widetilde{\Sigma}:=M\times\mathbbm{R}^{2m}\times S^{1}$. Suppose
$\mathcal{O}\subseteq M$ is open and $c_{\textrm{\emph{HZ}}}(\mathcal{O},M)<\infty$.
Choose $r_{0}>0$ such that
\begin{equation}
\left\lceil \pi r_{0}^{2}\right\rceil <\left\lceil c_{\textrm{\emph{HZ}}}(\mathcal{O},M)\right\rceil 
\end{equation}
and set 
\begin{equation}
r_{1}:=\sqrt{\tfrac{1}{\pi}c_{\textrm{\emph{HZ}}}(\mathcal{O},M)}+1
\end{equation}
Then there does not exist $\varphi\in\mathrm{Cont}_{0,c}(\widetilde{\Sigma},\alpha+\lambda_{\mathrm{std}})$
such that 
\begin{equation}
\varphi(\mathcal{O}\times B(r_{1})\times S^{1})\subset\mathcal{O}\times B(r_{0})\times S^{1}.
\end{equation}
\end{Thm}
The proof of Theorem \ref{prop:our HZ result_INTRO
} is given in Section \ref{sub:Applications}. See also Corollary \ref{Irie} for an application of Theorem \ref{prop:our HZ result_INTRO
}. 

Finally, following \cite[Section 1.7]{EliashbergKimPolterovich2006}
we investigate a rigidity phenomenon of positive contractible loops
of contactomorphisms. Suppose now that $(M,d\gamma)$ is the completion
of a Liouville domain $(M_1,d\gamma_1)$. Set $S:=\partial M_1$
and $\kappa:=\gamma|_{S}$, so that $(S,\kappa)$ is a contact manifold.
Abbreviate
\begin{equation}
M_r:=\begin{cases}
M_1\backslash(S\times(r,1), & 0<r<1,\\
M_1\cup_{S}(S\times[1,r], & r\geq1.
\end{cases}
\end{equation}
We can prove the following result, which roughly speaking says that
if $\widetilde{\mathrm{Cont}}_0(S,\kappa)$ is \emph{non-orderable}, so there exists a positive
contractible loop $\chi=\{\chi_{t}\}_{t\in S^{1}}\subset\mathrm{Cont}_{0}(S,\kappa)$
of contactomorphisms, then it is \emph{not }possible to homotope $\zeta$
through \emph{positive} loops to $\id_{S}$. In \cite[Theorem 1.11]{EliashbergKimPolterovich2006}
this was proved for $S=S^{2n-1}$.
\begin{Thm}
Set $c:=c_{M}(M_1)$ and assume that $c<\infty$. Suppose that  $\chi=\{\chi_{t}\}_{t\in S^{1}}\subset \mathrm{Cont}_{0}(S,\kappa)$
a positive contractible loop of contactomorphisms. Let $g_{t}:S\rightarrow(0,\infty)$
denote the contact Hamiltonian of $\chi$, and set 
\begin{equation}
\varepsilon:=\min_{(t,y)\in S^{1}\times S}g_{t}(y)>0.
\end{equation}
Then if $\{\chi_{s,t}\}_{0\leq s\leq1}$ is any homotopy of loops
of contactomorphisms such that $\chi_{1,t}=\chi_{t}$ and $\chi_{0,t}=\mbox{\emph{id}}_{S}$
with corresponding contact Hamiltonian $g_{s,t}:S\rightarrow\mathbbm{R}$
then there exists $(s,t,y)\in[0,1]\times S^{1}\times S$ such that
$g_{s,t}(y) \le -(1-\varepsilon)c$. 
\end{Thm}
\begin{proof}
This follows directly from Theorem \ref{thm:liouville flow} and the
material from \cite[Section 2.1]{EliashbergKimPolterovich2006}. Indeed,
suppose there exists $\delta>0$ such that $g_{s,t}(y)>-(1-\varepsilon)(c-\delta)$
for all $(s,t,y)\in[0,1]\times S^{1}\times S$. Set $a:=\min\{\varepsilon,\varepsilon c\}$.
Then as proved in \cite[Section 2.1]{EliashbergKimPolterovich2006}
for any $r<\frac{1}{c-\delta}$ it is possible to squeeze $M_r\times S^{1}$
into $M_{\frac{r}{1+ar}}\times S^{1}$. Fix $0<\lambda<\min\{a,\delta\}$
and take $r=\tfrac{1}{c-\lambda}$. Then 
\begin{align}
\overline{c}(M_r\times S^{1}) & =\left\lceil c_{M}(M_r)\right\rceil \nonumber \\
 & =\left\lceil rc_{M}(M_1)\right\rceil \nonumber \\
 & =\left\lceil \tfrac{c}{c-\lambda}\right\rceil =2.
\end{align}
But 
\begin{align}
\overline{c}(M_{\frac{r}{1+ar}})\times S^{1}) & =\left\lceil c_{M}(M_{\frac{r}{1+ar}})\right\rceil \nonumber \\
 & =\left\lceil \frac{cr}{1+ar}\right\rceil \nonumber \\
 & =\left\lceil \frac{c}{c-\lambda+a}\right\rceil =1.
\end{align}
This contradicts Theorem \ref{thm:liouville flow}.\end{proof}

\begin{acknowledgement*}
We are very grateful to Otto van Koert for explaining Example \ref{Ex:Otto} to us, and for many helpful discussions. We are also very grateful to Urs Fuchs, Leonid Polterovich, Daniel Rosen, Sheila Sandon and Frol Zapolsky for their helpful comments and useful remarks. Finally we thank Irida Altman for her help with Figure \ref{fig_irida}. PA is supported by the SFB 878 - Groups, Geometry and Actions, and WM is supported by an ETH Postdoctoral Fellowship. 
\end{acknowledgement*}

\section{\label{sec:Preliminaries}Preliminaries}

\subsection{\label{sub:Introductory-definitions}Introductory definitions}

Suppose $(\Sigma,\xi)$ is a connected closed coorientable contact
manifold. We denote by $\mathcal{P}\mathrm{Cont}_{0}(\Sigma,\xi)$ the
set of all smoothly parametrized paths $\{\varphi_{t}\}_{0\leq t\leq1}$
with $\varphi_{0}=\mathrm{id}_{\Sigma}$. We introduce an equivalence
relation $\sim$ on $\mathcal{P}\mathrm{Cont}_{0}(\Sigma,\xi)$ by saying
that two paths $\varphi$ and $\psi$ are equivalent if $\varphi_{1}=\psi_{1}$
and we can connect $\varphi$ and $\psi$ via a smooth family $\varphi^{s}=\{\varphi_{t}^{s}\}_{0\leq s,t\leq1}$
of paths such that $\varphi^{0}=\varphi$, $\varphi^{1}=\psi$ and
such that $\varphi_{1}^{s}$ is independent of $s$. The universal
cover $\widetilde{\mathrm{Cont}}_{0}(\Sigma,\xi)$ of $\mathrm{Cont}_{0}(\Sigma,\xi)$
is then $\mathcal{P}\mathrm{Cont}_{0}(\Sigma,\xi)/\sim$. We now give
the precise definition of a Liouville manifold, and what it means
for $\Sigma$ to be Liouville fillable.
\begin{Def}
\label{A-Liouville-manifold}A \emph{Liouville} \emph{domain }$(W_1,\lambda_1)$
is a compact exact symplectic manifold such that $\lambda_{1}|_{\partial W_1}$
is a positive contact form on $\partial W_1$. Equivalently the
vector field $Z_1$ on $W_1$ defined by $\iota_{Z_1}\lambda_1 =d\lambda_1$
should be transverse to $\partial W_1$ and point outwards. $Z_1$
is called the \emph{Liouville vector field}, and we denote by $l^t$
the flow of $Z_1$, which is defined for all $t\leq0$, and thus
induces an embedding $\partial W_1 \times(0,1]\rightarrow W_1$
defined by $(x,r)\mapsto l^{\log\, r}(x)$. Thus we can form the \emph{completion
}$(W,d\lambda)$ of $(W_1,\lambda_1)$ by attaching $\partial W_1\times[1,\infty)$
onto $\partial W_1$: 
\[
W:=W_1\cup_{\partial W_1}(\partial W_1\times[1,\infty)).
\]
We extend $\lambda_1$ to a 1-form $\lambda$ on $W$ by setting
$\lambda=r\lambda_1|_{\partial W_1}$ on $\partial W_1 \times[1,\infty)$.
Thus $d\lambda$ is a symplectic form on $W$. Similarly we extend
$Z_1$ to a vector field $Z$ on $W$ by setting $Z=r\partial_{r}$
on $\partial W_1\times[1,\infty)$. One calls $(W,d\lambda)$ a
\emph{Liouville manifold} - thus Liouville manifolds are exact non-compact
symplectic manifolds obtained by completing a Liouville domain. 

\label{def - liouville fillable}
We say that a closed connected coorientable
contact manifold $(\Sigma,\xi)$ is \emph{Liouville fillable }if there
exists a Liouville domain $(W_1,d\lambda_1)$ such that $\Sigma=\partial W_1$
and such that if $\alpha:=\lambda|_{\Sigma}$ then $\alpha$ is a
positive contact form on $\Sigma$ with $\ker\,\alpha=\xi$. By a
slight abuse of notation we will generally refer to the Liouville
manifold $(W,d\lambda)$ obtained from completing $(W_1,d\lambda_1)$
as ``the" filling of $\Sigma$. 
\end{Def}
The \emph{symplectization }$S\Sigma$ of a contact manifold $(\Sigma,\xi=\ker\,\alpha)$
is the symplectic manifold $\Sigma\times(0,\infty)$ equipped with
the symplectic form $d(r\alpha)$. If $\Sigma$ is Liouville fillable
with filling $(W,d\lambda)$ then one can embed $S\Sigma\hookrightarrow W$
by using the flow $l^{t}$ of the Liouville vector field $Z$ of $V$.
Next we recall how to lift a path $\varphi=\{\varphi_{t}\}_{0\leq t\leq1}$
to a symplectic isotopy $\Phi=\{\Phi_{t}\}_{0\leq t\leq1}$ on the
symplectization $S\Sigma$. Write $\varphi_{t}^{*}\alpha=\rho_{t}\varphi_{t}$.
Then define $\Phi_{t}:S\Sigma\rightarrow S\Sigma$ by 
\begin{equation}
\Phi_{t}(x,r):=\left(\varphi_{t}(x),\frac{r}{\rho_{t}(x)}\right).\label{eq:lift path to symplectization}
\end{equation}
The path $\Phi_{t}$ is Hamiltonian (in fact it preserves $\lambda$)
with Hamiltonian function 
\begin{equation}
H_{t}(x,r):=rh_{t}(x):S\Sigma\rightarrow\mathbbm{R}.\label{eq:Ht}
\end{equation}

We next define precisely what it means for a contact form $\alpha$
generating $\xi$ to be of Morse-Bott type.
\begin{Def}
\label{def morse bott type}A contact 1-form $\alpha\in\Omega^{1}(\Sigma)$
generating $\xi$ is said to be of \emph{Morse-Bott type }if for each
$T>0$, the set $P_{T}\subset\Sigma$ of points $x\in\Sigma$ satisfying
$\theta^{T}(x)=x$ is a closed submanifold of $\Sigma$, with the
property that $\mathrm{rank}\, d\alpha|_{P_{T}}$ is locally constant
and 
\begin{equation}
T_{x}P_{T}=\ker\left(D\theta^{T}(x)-\mathbbm{1}_{T_{x}\Sigma}\right) \qquad \mathrm{for\ all\ }x\in P_{T}.
\end{equation}
A Liouville fillable contact manifold $(\Sigma,\xi)$ is said to admit
a \emph{Morse-Bott Liouville filling }if there exists a filling $(W,d\lambda)$
such that $\alpha:=\lambda|_{\Sigma}$ is of Morse-Bott type. 
\end{Def}
Let us now recall the definition of a translated point of a contactomorphism.
This notion was introduced by Sandon in \cite{Sandon2012}. 
\begin{Def}
\label{def trans poin}
Let $(\Sigma,\xi)$ denote a closed connected
coorientable contact manifold, and fix a contact form $\alpha\in\Omega^{1}(\Sigma)$
generating $\xi$. Fix $\psi\in\mathrm{Cont}_{0}(\Sigma,\xi)$. We
can write $\psi^{*}\alpha=\rho\alpha$ for a smooth positive function
$\rho$ on $\Sigma$. A \emph{translated point }of $\psi$ is a
point $x\in\Sigma$ with the property that there exists $\eta\in\mathbbm{R}$
such that 
\begin{equation}
\psi(x)=\theta^{\eta}(x),\ \ \ \mathrm{and}\ \ \ \rho(x)=1.\label{eq:TP}
\end{equation}
We call $\eta$ the \emph{time-shift }of $x$. Note that if the leaf
$\{\theta^{t}(x)\}_{t\in\mathbbm{R}}$ is closed (which is always the
case when $\alpha$ is periodic) then the time-shift is not unique.
Indeed, if the leaf $\{\theta^{t}(x)\}_{t\in\mathbbm{R}}$ has period
$T>0$ then $\psi(x)=\theta^{\eta+\nu T}(x)$ for all $\nu\in\mathbbm{Z}$.

\end{Def}
Now let us define what it means for a translated point $x$ of an
element $[\varphi]\in\widetilde{\mathrm{Cont}}_{0}(\Sigma,\xi)$ to be
\emph{contractible }with respect to a Liouville filling $(W,d\lambda)$.
\begin{Def}
\label{def contractible translated point}Let $(W,d\lambda)$ denote
a Liouville filling of $(\Sigma,\xi)$, with $\alpha=\lambda|_{\Sigma}$.
Suppose $[\varphi]\in\widetilde{\mathrm{Cont}}_{0}(\Sigma,\xi)$ and $x$
is a translated point of $\varphi_1$ of time-shift $\eta$. We say
that the pair $(x,\eta)$ is a \emph{contractible }translated point
if the continuous loop $l:\mathbbm{R}/2\mathbbm{Z}\rightarrow\Sigma$
obtained from concatenating the path $\{\varphi_{t}(x)\}_{0\leq t\leq1}$
with the path $\{\theta^{-\eta t}(x)\}_{0\leq t\leq1}$ is contractible
in $W$. It is easy to see that this does not depend on path $\varphi=\{\varphi_{t}\}_{0\leq t\leq1}\in\mathcal{P}\mathrm{Cont}_{0}(\Sigma,\xi)$
representing $[\varphi]$. 
\end{Def}

For us, the usefulness of translated points stems from the fact that
the translated points of $\varphi$ are the generators of the Rabinowitz
Floer homology associated to $\varphi$, when the Rabinowitz Floer
homology is well defined; see Lemma \ref{Lemma:critical points} or
\cite{AlbersMerry2013a} for more information.

\subsection{The Rabinowitz action functional $\mathcal{A}_{\varphi}$ on $\Lambda(S\Sigma)\times\mathbbm{R}$}

Write $\Lambda(S\Sigma):=C_{\textrm{contr}}^{\infty}(S^{1},S\Sigma)$
for the component of the free loop space containing the contractible
loops.
\begin{Def}
Fix a path $\varphi\in\mathcal{P}\mathrm{Cont}_{0}(\Sigma,\xi)$ as
above, and let $H_{t}$ denote the Hamiltonian (\ref{eq:Ht}). We
define the \emph{perturbed Rabinowitz action functional} 
\begin{equation}
\mathcal{A}_{\varphi}:\Lambda(S\Sigma)\times\mathbbm{R}\rightarrow\mathbbm{R}
\end{equation}
by
\begin{equation}
\mathcal{A}_{\varphi}(u,\eta):=\int_{0}^{1}u^{*}\lambda-\eta\int_{0}^{1}\beta(t)(r(t)-1)dt-\int_{0}^{1}\dot{\chi}(t)H_{\chi(t)}(u(t))dt,\label{eq:Rab functional}
\end{equation}
where $\beta:\mathbbm{R}/\mathbbm{Z}\rightarrow\mathbbm{R}$ is a smooth
function with
\begin{equation}
\beta(t)=0\ \forall t\in[\tfrac{1}{2},1],\ \ \ \mathrm{and}\ \ \ \int_{0}^{1}\beta(t)dt=1,
\end{equation}
and $\chi:[0,1]\rightarrow[0,1]$ is a smooth monotone map with $\chi(\tfrac{1}{2})=0$,
$\chi(1)=1$, and $r(t)$ is the $\mathbbm{R}$-component of the map
$u:S^{1}\to S\Sigma=\Sigma\times\mathbbm{R}$. Denote by $\mathrm{Crit}(\mathcal{A}_{\varphi})$
the set of critical points of $\mathcal{A}_{\varphi}$, and denote
by $\mathrm{Spec}({\varphi}):=\mathcal{A}_{\varphi}(\mathrm{Crit}(\mathcal{A}_{\varphi}))$.\end{Def}
\begin{Rmk}
In this paper we define the Rabinowitz action functional only on the
contractible component of the free loop space, as all the applications
we have in mind here pertain only to the contractible component. Nevertheless,
it is possible to carry out all of our constructions on the full loop
space without any changes. This is because the symplectic form (on
$S\Sigma$ or on the Liouville filling $W$) is exact, and so there
are no ambiguities in the definition of $\mathcal{A}_{\varphi}$ on
non-contractible loops. 
\end{Rmk}
The following lemma explains why the perturbed Rabinowitz action functional
is useful in detecting translated points. It is a minor variant of
an argument originally due to the first author and Frauenfelder \cite[Proposition 2.4]{AlbersFrauenfelder2010c}.
For the convenience of the reader we recall the proof again here.
\begin{Lemma}
\label{Lemma:critical points}
\cite{AlbersMerry2013a}
A pair $(u,\eta)$ is a critical point of $\mathcal{A}_{\varphi}$ only if, writing
$u(t)=(x(t),r(t))\in\Sigma\times(0,\infty)$, $p:=x(\tfrac{1}{2})$
is a translated point of $\varphi$, with time-shift $-\eta$. Conversely every such pair $( p, \eta)$ gives rise to a unique critical point of $ \mathcal{A}_{ \varphi}$. Moreover
if $(u,\eta)$ is a critical point of $\mathcal{A}_{\varphi}$ then
\begin{equation}
\mathcal{A}_{\varphi}(u,\eta)=\eta.\label{eq:value on crit}
\end{equation}
If $\varphi$ is an exact path of contactomorphisms then $r(t)\equiv1$
for every critical point $(u=(x,r),\eta)$. 
\end{Lemma}
\begin{proof}
Denote by $\Phi_{t}:S\Sigma\rightarrow S\Sigma$ the symplectic isotopy
(\ref{eq:lift path to symplectization}). A pair $(u,\eta)$ with
$u=(x,r):S^{1}\rightarrow\Sigma\times(0,\infty)$ belongs to $\mathrm{Crit}(\mathcal{A}_{\varphi})$
if and only if 
\begin{equation}
\begin{cases}
\dot{u}(t)=\eta\beta(t)R(x(t))+\dot{\chi}(t)X_{H_{\chi(t)}}(u(t)),\\
\int_{0}^{1}\beta(t)(r(t)-1)dt=0.
\end{cases}\label{eq:critical points}
\end{equation}
Thus for $t\in[0,\tfrac{1}{2}]$, we have $r(t)=1$ and $\dot{x}(t)=-\eta R(x(t))$,
and $x(1)=\Phi_{1}(x(\tfrac{1}{2}))$. Suppose $(u,\eta)\in\mathrm{Crit}(\mathcal{A}_{\varphi})$.
Thus $u(\tfrac{1}{2})=(\theta^{-\eta}(x(0)),1)$. Next, for $t\in[\tfrac{1}{2},1]$
we have $\dot{u}(t)=\dot{\chi}(t)X_{H_{\chi(t)}}(u(t))$. In particular,
$\varphi(x(\tfrac{1}{2}))=\theta^{-\eta}(x(\tfrac{1}{2}))$, and thus
$x(\tfrac{1}{2})$ is a translated point of $\varphi$. Moreover the
time shift of $x$ is given by $-\eta$ mod 1.

Next, we note that
\begin{equation}
\lambda(X_{H}(x,r))=dH(x,r)\left(r\frac{\partial}{\partial r}\right)=H(x,r),
\end{equation}
and hence
\begin{align}
\mathcal{A}_{\varphi}(u,\eta) & =\int_{0}^{\tfrac{1}{2}}(r\alpha)(\eta\beta(u)R(x))dt+\int_{\tfrac{1}{2}}^{1}\left[\lambda(\dot{\chi}X_{H_{\chi}}(u))-\dot{\chi}H_{\chi}(u)\right]dt\nonumber \\
 & =\eta+0.
\end{align}
Finally if $\varphi_{t}$ is exact for each $t$ then the path $\Phi_{t}$
of symplectomorphisms defined in (\ref{eq:lift path to symplectization})
is simply given by $\Phi_{t}(x,r)=(\varphi_{t}(x),r)$, from which
the last statement immediately follows.\end{proof}
\begin{Rmk}
\label{the ubiquitous minus sign}
We emphasize again that if $(u=(x,r),\eta)$
is a critical point of $\mathcal{A}_{\varphi}$ then the time-shift
of the translated point $x(\tfrac{1}{2})$ is the \emph{negative}
of the action value. This explains the Reeb flow will turn out to
have a \emph{negative }spectral number (cf. part (1) of Theorem \ref{Theorem-A}). 

We point out that there is a distinguished Morse-Bott component of  $ \mathrm{Crit}(\mathcal{A}_{ \mathrm{id}_{ \Sigma}})$ diffeomorphic to $ \Sigma$  corresponding to critical points $((x,1), 0)$ for $x \in \Sigma$.
\end{Rmk}

We now define what if means for $\varphi$ to be non-degenerate. In
the periodic case we also introduce the notion of being non-resonant. 
\begin{Def}
\label{def: nondegen}A path $\varphi$ is \emph{non-degenerate }if
$\mathcal{A}_{\varphi}:\Lambda(S\Sigma)\times\mathbbm{R}\rightarrow\mathbbm{R}$
is a Morse-Bott function. In the periodic case we say that $\varphi$
is \emph{non-resonant }if $\mathrm{Spec}({\varphi})\cap\mathbbm{Z}=\emptyset$.
\end{Def}
\begin{Rmk}
\label{Rmk: spec only depends on the end}
The identity $\mathrm{id}_{\Sigma}$ is non-degenerate if and only if
$\alpha$ is of Morse-Bott type (see \cite[Appendix B]{CieliebakFrauenfelder2009}).
In \cite{AlbersMerry2013a} we explained why a generic path $\varphi$
is non-degenerate (actually a stronger result is true: for generic
$\varphi$ the functional $\mathcal{A}_{\varphi}$ is even Morse).
It is also easy to see that a generic $\varphi$ is non-resonant. 
Moreover $ \mathrm{Spec}({ \varphi})$ is always a nowhere dense subset of $ \R$ (even in the degenerate case), cf. \cite[Lemma 3.8]{Schwarz2000}. Finally we note that $\mathrm{Spec}({\varphi})$
depends only on the terminal map $\varphi_{1}$.
\end{Rmk}
The next lemma explains why we pay particular attention to periodic
contact manifolds. It will prove crucial in the construction of the
contact capacity (cf. Section \ref{sec:The-periodic-case - capacities},
in particular Proposition \ref{lem the use of periodicity}). Fix
$\varphi\in\widetilde{\mathrm{Cont}}_{0}(\Sigma,\xi)$ and fix a contactomorphism
$\psi\in\mathrm{Cont}_{0}(\Sigma,\xi)$.
\begin{Lemma}\label{lem:integral_critical_points_are_conjugation_invariant}
Assume
$\alpha$ is periodic. If $(u=(x,r),\eta)\in\mathrm{Crit}(\mathcal{A}_{\varphi})$
with $\eta\in\mathbbm{Z}$ then there exists a critical point $(u_{1}=(x_{1},r_{1}),\eta)$
of $\mathcal{A}_{\psi\varphi\psi^{-1}}$ with $x_{1}(\tfrac{1}{2})=\psi(x(\tfrac{1}{2}))$.
In particular, 
\begin{equation}
\mathrm{Spec}({\varphi})\cap\mathbbm{Z}=\emptyset\ \ \ \Leftrightarrow\ \ \ \mathrm{Spec}({\psi\varphi\psi^{-1}})\cap\mathbbm{Z}=\emptyset.\label{eq:key eq}
\end{equation}
Moreover $(u,\eta)$ is non-degenerate if and only if $(u_{1},\eta)$
is non-degenerate. \end{Lemma}
\begin{proof}
If $(u,\eta)\in\mathrm{Crit}(\mathcal{A}_{\varphi})$ with $\eta\in\mathbbm{Z}$
then since $\theta^{t}$ is 1-periodic, this means that if we write
$u(t)=(x(t),r(t))$ then $x(\tfrac{1}{2})$ is a\emph{ }fixed point
of $\varphi$. Thus $\psi(x(\tfrac{1}{2}))$ is a fixed point of $\psi\varphi\psi^{-1}$.
Thus by Lemma \ref{Lemma:critical points} for each $\nu\in\mathbbm{Z}$
there exists a critical point $(u_{\nu}=(x_{\nu},r_{\nu}),\nu)$ of
$\mathcal{A}_{\psi\varphi\psi^{-1}}$ with $x_{\nu}(\tfrac{1}{2})=\psi(x_{\nu}(\tfrac{1}{2}))$.
In particular, this is true for $\nu=\eta$.
The final statement follows from the fact that the linearised equation is also conjugation invariant.  
\end{proof}

\subsection{\label{sec:Rabinowitz-Floer-homology}Rabinowitz Floer homology}

Let us now assume that $(\Sigma,\xi)$ is Liouville fillable with
a Morse-Bott Liouville filling $(W,d\lambda)$. We would like to extend
$\mathcal{A}_{\varphi}$ to a functional defined on all of $\Lambda(W)\times\mathbbm{R}$,
where $\Lambda(W):=C_{\textrm{contr}}^{\infty}(S^{1},W)$ as before.
In order to do this we must extend the function $(x,r)\mapsto r-1$
and the Hamiltonian $H_{t}$ to functions defined on all of $W$.
At the same time, it is convenient to truncate them. As in \cite{AlbersFrauenfelder2010c},
we proceed as follows. Define $m:W\rightarrow\mathbbm{R}$ so that
\begin{equation}
m(x,r):=r-1\ \ \ \mathrm{on\ }\Sigma\times(\tfrac{1}{2},\infty),
\end{equation}
\begin{equation}
\frac{\partial m}{\partial r}(x,r)\geq0\ \ \ \mathrm{for\ all\ }(x,r)\in S\Sigma,
\end{equation}
\begin{equation}
m|_{W\backslash S\Sigma}:=-\tfrac{3}{4}.
\end{equation}
Next, for $\kappa>0$ let $\varepsilon_{\kappa}\in C^{\infty}([0,\infty),[0,1])$
denote a smooth function such that 
\begin{equation}
\varepsilon_{\kappa}(r)=\begin{cases}
1, & r\in[e^{-\kappa},e^{\kappa}],\\
0, & r\in[0,e^{-2\kappa}]\cup[e^{\kappa}+1,\infty),
\end{cases}
\end{equation}
and such that 
\begin{equation}
0\leq\varepsilon_{\kappa}'(r)\leq2e^{2\kappa}\ \ \ \mathrm{for\ }r\in[e^{-2\kappa},e^{-\kappa}],
\end{equation}
\begin{equation}
-2\leq\varepsilon_{\kappa}'(r)\leq0\ \ \ \mathrm{for\ }r\in[e^{\kappa},e^{\kappa}+1].
\end{equation}
Then define $H_{t}^{\kappa}:W\rightarrow\mathbbm{R}$ by setting $H_{t}^{\kappa}|_{W\backslash S\Sigma}=-\tfrac{3}{4}$
and
\begin{equation}
H_{t}^{\kappa}(x,r):=\varepsilon_{\kappa}(r)H_{t}(x,r)\ \ \ \mathrm{for\ }(x,r)\in S\Sigma.\label{eq:truncating H}
\end{equation}
We denote by $\mathcal{A}_{\varphi}^{\kappa}:\Lambda(W)\times\mathbbm{R}\rightarrow\mathbbm{R}$
the Rabinowitz action functional defined using the Hamiltonians $m$
and $H_{t}^{\kappa}$: 
\begin{equation}
\mathcal{A}_{\varphi}^{\kappa}(u,\eta):=\int_{0}^{1}u^{*}\lambda-\eta\int_{0}^{1}\beta(t)m(u(t))dt-\int_{0}^{1}\dot{\chi}(t)H_{\chi(t)}^{\kappa}(u(t))dt.
\end{equation}

\begin{Def}
\label{hofer norm def}Now let us recall the definition of the oscillation
`norm' on $\widetilde{\mathrm{Cont}}_{0}(\Sigma,\xi)$. Firstly, suppose
$\{\varphi_{t}\}_{0\leq t\leq1}\in\mathcal{P}\mathrm{Cont}_{0}(\Sigma,\xi)$.
Let $h_{t}:\Sigma\rightarrow\mathbbm{R}$ denote the contact Hamiltonian.
The \emph{oscillation norm} $\left\Vert h\right\Vert _{\textrm{osc}}$
is defined by 
\begin{equation}
\left\Vert h\right\Vert _{\textrm{osc}}:=\left\Vert h\right\Vert _{+}+\left\Vert h\right\Vert _{-},\label{eq:osc 1}
\end{equation}
\begin{equation}
\left\Vert h\right\Vert _{+}:=\int_{0}^{1}\max_{x\in\Sigma}h_{t}(x)dt,\ \ \ \left\Vert h\right\Vert _{-}:=-\int_{0}^{1}\min_{x\in\Sigma}h_{t}(x)dt.\label{eq:osc 2}
\end{equation}
We then define the oscillation norm\emph{}%
\footnote{Note that this is not a norm; the same is true of the quantity $\left\Vert \varphi\right\Vert _{\kappa}$
defined below. %
} $\left\Vert \varphi\right\Vert $ and its positive and negative parts
$\left\Vert \varphi\right\Vert _{\pm}$ for $\varphi\in\widetilde{\mathrm{Cont}}_{0}(\Sigma,\xi)$
by taking the infimum of the oscillation norms $\left\Vert h\right\Vert _{\textrm{osc}}$
{[}resp. $\left\Vert h\right\Vert _{\pm}${]} over all possible paths
$\{\varphi_{t}\}_{0\leq t\leq1}$ representing $\varphi$ (with corresponding
contact Hamiltonians $h_{t}$).\end{Def}
\begin{Rmk}
Denote by $\Phi^{\kappa}\in\mathrm{Ham}_{c}(W,d\lambda)$ the Hamiltonian
diffeomorphism generated by $H_{t}^{\kappa}$. The\textbf{\emph{ }}\emph{Hofer
norm} $\left\Vert \Phi^{\kappa}\right\Vert $ of $H_{t}^{\kappa}$
is related to $\left\Vert \varphi\right\Vert $ by 
\begin{equation}
\left\Vert \Phi^{\kappa}\right\Vert \leq e^{\kappa}\left\Vert \varphi\right\Vert .
\end{equation}
Here by definition the Hofer norm $\left\Vert \Phi^{\kappa}\right\Vert $
is the infimum of the oscillation norms of all possible Hamiltonians
generating $\Phi^{\kappa}$. One such Hamiltonian is $H_{t}^{\kappa}$,
and it is clear that $\left\Vert H^{\kappa}\right\Vert _{\textrm{osc}}\leq e^{\kappa}\left\Vert h\right\Vert _{\textrm{osc}}$. \end{Rmk}
\begin{Def}
Suppose $\varphi\in\mathcal{P}\mathrm{Cont}_{0}(\Sigma,\xi)$. Let $\rho_{t}:\Sigma\rightarrow(0,\infty)$
is defined by $\varphi_{t}^{*}\alpha=\rho_{t}\alpha$. Define a constant
$\kappa(\varphi)>0$ by 
\begin{equation}
\label{eq:c phi}
\kappa(\varphi):=\max_{t\in[0,1]}\left|\int_{0}^{t}\max_{x\in\Sigma}\frac{\dot{\rho}_{\tau}(x)}{\rho_{\tau}(x)^{2}}d\tau\right|
\end{equation}
Note that if $\varphi$ is exact then $\kappa(\varphi)=0$. 
\end{Def}
In \cite[Proposition 2.5]{AlbersMerry2013a} we proved:
\begin{Lemma}
\label{prop:linfinity-1}If $\kappa>\kappa(\varphi)$ then if $(u,\eta)\in\mathrm{Crit}(\mathcal{A}_{\varphi}^{\kappa})$
then $u(S^{1})\subseteq S\Sigma$, and moreover if we write $u(t)=(x(t),r(t))$
then $r(S^{1})\subseteq(e^{-\kappa/2},e^{\kappa/2})$. 
\end{Lemma}

If $\varphi$ is non-degenerate in the sense of Definition \ref{def: nondegen},
then as explained in \cite{AlbersMerry2013a}, Lemma \ref{prop:linfinity-1}
allows to define for $a,b\in[-\infty,\infty]\backslash\mathrm{Spec}({\varphi})$
the \emph{Rabinowitz Floer homology}
\begin{equation}
\RFH_{*}^{(a,b)}(\mathcal{A}_{\varphi},W).
\end{equation}
This is a semi-infinite dimensional Morse theory associated to the
functional $\mathcal{A}_{\varphi}^{\kappa}$ (for some $\kappa>\kappa_{0}(\varphi)$),
and we sketch the definition here and refer to e.g. \cite{AlbersFrauenfelder2010c,AlbersMerry2013a}
for more information.

Let $\mathcal{J}_{\textrm{conv}}(S\Sigma)$ denote the set of time
dependent almost complex structures $J=\{J_{t}\}_{t\in S^{1}}$ on
$S\Sigma$ that are $d(r\alpha)$-compatible and that are \emph{convex}.
Here we use the sign convention that $J$ is $d\lambda$-compatible
if $d(r\alpha)(J\cdot,\cdot)$ defines a Riemannian metric, and to
say that $J$ is convex is to ask that there exists $S_{0}>0$ such
that 
\begin{equation}
dr\circ J_{t}=d(r\alpha)\ \ \ \mathrm{on}\ \ \Sigma\times[S_{0},\infty)\label{eq:convex type}
\end{equation}
(in particular, $J$ is independent of $t$ on $\Sigma\times[S_{0},\infty)$).
We denote by $\mathcal{J}_{\textrm{conv}}(W)$ the set of time-dependent
almost complex structures $J=\{J_{t}\}_{t\in S^{1}}$ with the property
that $J|_{S\Sigma}\in\mathcal{J}_{\textrm{conv}}(S\Sigma)$. 

Given $J\in\mathcal{J}_{\textrm{conv}}(W)$ we can define an $L^{2}$-inner
product $\left\langle \left\langle \cdot,\cdot\right\rangle \right\rangle _{J}$
on $\Lambda(W)\times\mathbbm{R}$: for $(u,\eta)\in\Lambda(W)\times\mathbbm{R}$,
$\zeta,\zeta'\in\Gamma(u^{*}TW)$ and $b,b'\in\mathbbm{R}$, set 
\begin{equation}
\left\langle \left\langle (\zeta,b),(\zeta',b')\right\rangle \right\rangle _{J}:=\int_{0}^{1}d\lambda(J_{t}\zeta,\zeta')dt+bb'.\label{eq:l2 norm}
\end{equation}
We denote by $\nabla_{J}\mathcal{A}_{\varphi}^{\kappa}$ the gradient
of $\mathcal{A}_{\varphi}^{\kappa}$ with respect to $\left\langle \left\langle \cdot,\cdot\right\rangle \right\rangle _{J}$. 
\begin{Rmk}
In this paper all sign conventions are the same as in \cite{AbbondandoloSchwarz2009,AlbersMerry2013a}.
\end{Rmk}
Assume that $\varphi$ is non-degenerate and fix $\kappa>\kappa_{0}(\varphi)$
and $J\in\mathcal{J}_{\textrm{conv}}(W)$. By assumption $\mathcal{A}_{\varphi}^{\kappa}$
is a Morse-Bott function.\emph{ }Pick a Morse function $g:\mathrm{Crit}(\mathcal{A}_{\varphi}^{\kappa})\rightarrow\mathbbm{R}$,
and choose a Riemannian metric $\varrho$ on $\mathrm{Crit}(\mathcal{A}_{\varphi}^{\kappa})$
such that the negative gradient flow of $\nabla_{\varrho}g$ is \emph{Morse-Smale}.
Given two critical points $w^{-},w^{+}\in\mathrm{Crit}(g)$, with $w^{\pm}=(u^{\pm},\eta^{\pm})$,
we denote by $\mathcal{M}_{w^{-},w^{+}}(\mathcal{A}_{\varphi}^{\kappa},g,J,\varrho)$
the moduli space of \emph{gradient flow lines with cascades} of $-\nabla_{J}\mathcal{A}_{\varphi}^{\kappa}$
and $-\nabla_{\varrho}g$ running from $w^{-}$ to $w^{+}$. See \cite[Appendix A]{Frauenfelder2004}
or \cite[Appendix A]{CieliebakFrauenfelder2009} for the precise definition.

Introduce a grading on $\mathrm{Crit}(g)$ by setting 
\begin{equation}
\mu(u,\eta):=\begin{cases}
\mu_{\textrm{CZ}}(u)-\frac{1}{2}\dim_{(u,\eta)}\mathrm{Crit}(\mathcal{A}_{\varphi}^{\kappa})+\mathrm{ind}_{g}(u,\eta), & \eta>0,\\
\mu_{\textrm{CZ}}(u)-\frac{1}{2}\dim_{(u,\eta)}\mathrm{Crit}(\mathcal{A}_{\varphi}^{\kappa})+\mathrm{ind}_{g}(u,\eta)+1, & \eta<0,\\
1-n+\mathrm{ind}_{g}(u,\eta), & \eta=0.
\end{cases}\label{eq:grading}
\end{equation}
Here $\mu_{\textrm{CZ}}(u)$ denotes the Conley-Zehnder index of the
loop $t\mapsto u(t/\eta)$ and $\dim_{(u,\eta)}\mathrm{Crit}(\mathcal{A}_{\varphi}^{\kappa})$
denotes the local dimension of $\mathrm{Crit}(\mathcal{A}_{\varphi}^{\kappa})$
at $(u,\eta)$. Actually in most cases of interest in this paper,
we may assume that $\mathcal{A}_{\varphi}$ is actually Morse. In
this case the Morse function $g$ is taken to be identically zero,
and (\ref{eq:grading}) continues to hold. 
\begin{Rmk}
\label{Our-grading-convention} Our normalization convention for the
Conley-Zehnder index is that if $H$ is a $C^{2}$-small Morse function
on $W$ and $x$ is a critical point of $W$ then 
\begin{equation}
\mu_{\textrm{CZ}}(x)=n-\mathrm{ind}_{H}(x),
\end{equation}
where $\mathrm{ind}_{H}(x)$ denotes the Morse index of $x$. 
\end{Rmk}
Given $-\infty<a<b<\infty$ denote by $\mathrm{RFC}_{*}^{(a,b)}(\mathcal{A}_{\varphi}^{\kappa},g):=\mathrm{Crit}_{*}^{(a,b)}(g)\otimes\mathbbm{Z}_{2}$,
where $\mathrm{Crit}_{*}^{(a,b)}(g)$ denotes the set of critical points
$w$ of $g$ with $a<\mathcal{A}_{\varphi}(w)<b$. We only do this
when $a,b\notin\mathrm{Spec}({\varphi})$, even if this
is not explicitly stated. Generically the moduli spaces $\mathcal{M}_{w^{-},w^{+}}(\mathcal{A}_{\varphi},g,J,\varrho)$
carry the structure of finite dimensional smooth manifolds, whose
components of dimension zero are compact. One defines a boundary operator
$\partial$ on $\mathrm{RFC}_{*}^{(a,b)}(\mathcal{A}_{\varphi}^{\kappa},g)$
by counting the elements of the zero-dimensional parts of the moduli
spaces $\mathcal{M}_{w^{-},w^{+}}(\mathcal{A}_{\varphi},g,J,\varrho)$. 

The homology $\RFH_{*}^{(a,b)}(\mathcal{A}_{\varphi},W):=\mathrm{H}_{*}(\mathrm{RFC}_{*}^{(a,b)}(\mathcal{A}_{\varphi}^{\kappa},g),\partial)$
does not depend on any of the auxiliary choices we made. We emphasize
though that $\RFH_{*}^{(a,b)}(\mathcal{A}_{\varphi},W)$ depends
on the choice of filling $(W,d\lambda)$. Finally we define 
\begin{equation}
\RFH_{*}^{b}(\mathcal{A}_{\varphi},W):=\underset{a\downarrow-\infty}{\underrightarrow{\lim}}\RFH_{*}^{(a,b)}(\mathcal{A}_{\varphi},W),
\end{equation}
\begin{equation}
\RFH_{*}^{(a,\infty)}(\mathcal{A}_{\varphi},W):=\underset{b\uparrow\infty}{\underleftarrow{\lim}}\RFH_{*}^{(a,b)}(\mathcal{A}_{\varphi},W),
\end{equation}
\begin{equation}
\RFH_{*}(\mathcal{A}_{\varphi},W):=\underset{a\downarrow-\infty}{\underrightarrow{\lim}}\underset{b\uparrow\infty}{\underleftarrow{\lim}}\RFH_{*}^{(a,b)}(\mathcal{A}_{\varphi},W)\label{eq:inverse and direct limit}
\end{equation}
(the order of the limits in (\ref{eq:inverse and direct limit}) matters).
As pointed out by Ritter \cite{Ritter2013}, it follows from work
of Cieliebak-Frauenfelder-Oancea \cite{CieliebakFrauenfelderOancea2010}
that the Rabinowitz Floer homology $\RFH_{*}(\Sigma,W)$ vanishes
if and only if the symplectic homology $\mathrm{SH}_{*}(W)$ vanishes.\newline 

We briefly summarize now the key properties that we will need about
the Rabinowitz Floer homology $\RFH_{*}^{(a,b)}(\mathcal{A}_{\varphi},W)$:
\begin{enumerate}
\item \label{enu:The-Rabinowitz-Floer1}The Rabinowitz Floer homology is
independent of $\varphi$ in the following strong sense. There is
a universal object $\RFH_{*}(\Sigma,W)$ (which may be thought
as corresponding to the case $\varphi=\mathrm{id}_{\Sigma}$) together
with canonical isomorphisms 
\begin{equation}
\zeta_{\varphi}:\RFH_{*}(\Sigma,W)\rightarrow\RFH_{*}(\mathcal{A}_{\varphi},W).
\end{equation}
Given two paths $\varphi$ and $\psi$, there is a map $\zeta_{\varphi,\psi}:\RFH_{*}(\mathcal{A}_{\varphi},W)\rightarrow\RFH_{*}(\mathcal{A}_{\psi},W)$
with the property that 
\begin{equation}
\zeta_{\psi}=\zeta_{\varphi,\psi}\circ\zeta_{\varphi}.\label{eq:continuation homomorphisms}
\end{equation}
In particular, if $Z \in \RFH_*(\Sigma,W)$ is a non-zero class then $ \RFH_*( \mathcal{A}_{ \varphi},W)$ contains a non-zero class $ Z_{ \varphi}$ defined by 
\begin{equation}
\label{eq:Z_varphi}
\zeta_{\varphi,\psi}\big(Z_\varphi\big)=Z_\psi\quad\text{and}\quad Z_{ \mathrm{id}_{ \Sigma}} =Z \in\RFH_*(\Sigma,W).
\end{equation}

\item \label{enu:The-Rabinowitz-Floer2}If $a\leq b\leq\infty$ there is a natural map 
\begin{equation}
j_{\varphi}^{a,b}:\RFH_{*}^{a}(\mathcal{A}_{\varphi},W)\rightarrow\RFH_{*}^{b}(\mathcal{A}_{\varphi},W).
\end{equation}
 Similarly there is a natural map 
\begin{equation}
\label{proj_map}
p_{\varphi}^{a,b}:\RFH_{*}^{b}(\mathcal{A}_{\varphi},W)\rightarrow\RFH_{*}^{(a,b)}(\mathcal{A}_{\varphi},W).
\end{equation}
If $b=\infty$ we abbreviate $j_{\varphi}^{a,\infty}=j_{\varphi}^{a}$,
and we write $j^{a}$ for the map $\RFH_{*}^{a}(\Sigma,W)\rightarrow\RFH_{*}(\Sigma,W)$,
with similar conventions for the maps $p_{\varphi}^{a,b}$. If $\mathrm{Spec}({\varphi})\cap[a,b]=\emptyset$
then the map $j_{\varphi}^{a,b}:\RFH_{*}^{a}(\mathcal{A}_{\varphi},W)\rightarrow\RFH_{*}^{b}(\mathcal{A}_{\varphi},W)$
is an isomorphism and $p_{\varphi}^{a,b}:\RFH_{*}^{b}(\mathcal{A}_{\varphi},W)\rightarrow\RFH_{*}^{(a,b)}(\mathcal{A}_{\varphi},W)$
is the zero map (as $\RFH_{*}^{(a,b)}(\mathcal{A}_{\varphi},W)=0$). 
\item Moreover there is a filtered version of (\ref{eq:continuation homomorphisms}),
which gives the existence of a maps 
\begin{equation}
\label{eq:fil cont hom}
\zeta_{\varphi,\psi}^{a}:\RFH_{*}^{a}(\mathcal{A}_{\varphi},W)\rightarrow\RFH_{*}^{a+K(\varphi,\psi)}(\mathcal{A}_{\psi},W)
\end{equation}
for some constant $K(\varphi,\psi)\ge0$. The maps \eqref{eq:fil cont hom}
are a special case of \cite[Lemma 2.7]{AlbersFrauenfelder2010c}.
It will be important however to note that if the paths $\varphi,\psi$
have contact Hamiltonians $h_{t}$ and $k_{t}$ then then the constant
$K(\varphi,\psi)$ satisfies 
\begin{equation}
K(\varphi,\psi)\leq e^{\max \{ \kappa(\varphi), \kappa( \psi)  \} }\max\big\{\left\| h-k\right\| _{+},0\big\},\label{eq:the constant K}
\end{equation}
where we are using the notation from (\ref{eq:osc 1})-(\ref{eq:osc 2}).
Finally one has for all $Z\in\RFH_{*}^{a}(\mathcal{A}_{\varphi},W)$
that
\begin{equation}
\zeta_{\varphi,\psi}\left(j_{\varphi}^{a}(Z)\right)=j_{\psi}^{a+K(\varphi,\psi)}\left(\zeta_{\varphi,\psi}^{a}(Z)\right).\label{eq:commutativity}
\end{equation}
\item 
We recall from Remark \ref{the ubiquitous minus sign} that $\mathrm{Crit}( \mathcal{A}_{ \mathrm{id}_{ \Sigma}})$ contains  $ \Sigma$ as a Morse-Bott component via the constants. 
For $ \varepsilon > 0$ smaller than the smallest period of a contractible Reeb orbit, this gives rise to a canonical isomorphism
\begin{equation}
\label{eq:singular_homology_canonical}
  \RFH^{( - \varepsilon, \varepsilon)}_*( \Sigma, W) \cong \mathrm{H}_{ * + n-1}( \Sigma; \Z_2).
\end{equation}
\end{enumerate}

Even though it is more or less standard, the estimate \eqref{eq:the constant K} is extremely important in all that follows, and hence we prove it here. 
To define the continuation homomorphism $\zeta_{\varphi,\psi}$ we denote by $H_t=rh_t$ and $K_t=rk_t$ the Hamiltonian functions of $\varphi$ and $\psi$, respectively, and choose a linear homotopy
\begin{equation}
L^s_t:=\nu(s)H_t +(1-\nu(s))K_t
\end{equation}
for a smooth function $\nu:\R\to[0,1]$ with $\nu(s)=1$ for $s\leq-1$, $\nu(s)=0$ for $s\geq1$ and $\nu'(s)\leq0$. We define the ($s$-dependent) action functional $ \mathcal{A}_s $ as in \eqref{eq:Rab functional}:
\begin{equation}
\mathcal{A}_{s}(u,\eta)=\int_{0}^{1}u^{*}\lambda-\eta\int_{0}^{1}\beta(t)m(u(t))dt-\int_{0}^{1}\dot{\chi}(t)\varepsilon_{ \kappa}(r)L_{\chi(t)}^s(u(t))dt.
\end{equation}
where $ \varphi_s$ has corresponding Hamiltonian function $L^s_t$. 
Then counting solutions of 
\begin{equation}
\partial_{s}(u,\eta)+\nabla_{J}\mathcal{A}_{s}(u,\eta)=0
\end{equation}
with $(u_-,\eta_-):=\big(u(-\infty),\eta(-\infty)\big)\in\mathrm{Crit}(\mathcal{A}_{\varphi})$ and $(u_+,\eta_+):=\big(u(+\infty),\eta(+\infty)\big)\in\mathrm{Crit}(\mathcal{A}_{\psi})$
defines the continuation homomorphism. We recall that $\kappa>\max\{\kappa(\varphi),\kappa(\psi)\}$ and estimate
\begin{align}
0&\leq\mathbbm{E}_{J}(u,\eta)\\
&=\int_{-\infty}^{\infty}\int_{0}^{1}|\partial_{s}(u,\eta)|_{J}^{2}dtds\\
&=-\int_{-\infty}^{\infty}\int_{0}^{1}\left\langle \left\langle \nabla \mathcal{A}_{s}(u,\eta),\partial_{s}(u,\eta)\right\rangle \right\rangle _{J}dtds\\
&=-\int_{-\infty}^{\infty}\int_{0}^{1}\frac{d}{ds} \mathcal{A}_s(u,\eta)dtds+\int_{-\infty}^{\infty}\int_{0}^{1}\frac{\partial\mathcal{A}_s}{\partial s}(u,\eta)dtds\\
&=\mathcal{A}_{\varphi}(u_-,\eta_-)-\mathcal{A}_{\psi}(u_+,\eta_+)-\int_{-\infty}^{\infty}\int_{0}^{1}\dot{\chi}(t)\varepsilon_{\kappa}(r(t))\frac{\partial L^s_{\chi(t)}}{\partial s}(u(t))dtds\\
&=\mathcal{A}_{\varphi}(u_-,\eta_-)-\mathcal{A}_{\psi}(u_+,\eta_+)\\
&\qquad\qquad-\int_{-\infty}^{\infty}\int_{0}^{1}\nu'(s)\varepsilon_{\kappa}(r(t))\dot{\chi}(t)\Big(H_{\chi(t)}(u(t))-K_{\chi(t)}(u(t))\Big)dtds\\
&=\mathcal{A}_{\varphi}(u_-,\eta_-)-\mathcal{A}_{\psi}(u_+,\eta_+)\\
&\qquad\qquad-\int_{-\infty}^{\infty}\int_{0}^{1}\underbrace{\nu'(s)}_{\leq 0}\underbrace{\varepsilon_{\kappa}(r(t))r(t)}_{0\leq\cdot\leq e^\kappa}\underbrace{\dot{\chi}(t)}_{\geq0}\Big(h_{\chi(t)}(u(t))-k_{\chi(t)}(u(t))\Big)dtds\\
&\leq\mathcal{A}_{\varphi}(u_-,\eta_-)-\mathcal{A}_{\psi}(u_+,\eta_+)\\
&\qquad\qquad-\int_{-\infty}^{\infty}\int_{0}^{1}\underbrace{\nu'(s)}_{\leq 0}\underbrace{\varepsilon_{\kappa}(r(t))r(t)}_{0\leq\cdot\leq e^\kappa}\underbrace{\dot{\chi}(t)}_{\geq0}\max_{x\in\Sigma}\Big(h_{\chi(t)}(x)-k_{\chi(t)}(x)\Big)dtds\\
&\leq\mathcal{A}_{\varphi}(u_-,\eta_-)-\mathcal{A}_{\psi}(u_+,\eta_+)\\
&\qquad\qquad-\int_{-\infty}^{\infty}\int_{0}^{1}\underbrace{\nu'(s)}_{\leq 0}\underbrace{\varepsilon_{\kappa}(r(t))r(t)}_{0\leq\cdot\leq e^\kappa}\underbrace{\dot{\chi}(t)}_{\geq0}\max\Big\{\max_{x\in\Sigma}\Big(h_{\chi(t)}(x)-k_{\chi(t)}(x)\Big),0\Big\}dtds\\
&\leq\mathcal{A}_{\varphi}(u_-,\eta_-)-\mathcal{A}_{\psi}(u_+,\eta_+)\\
&\qquad\qquad-e^\kappa\int_{-\infty}^{\infty}\int_{0}^{1}\nu'(s)\dot{\chi}(t)\max\Big\{\max_{x\in\Sigma}\Big(h_{\chi(t)}(x)-k_{\chi(t)}(x)\Big),0\Big\}dtds\\
&=\mathcal{A}_{\varphi}(u_-,\eta_-)-\mathcal{A}_{\psi}(u_+,\eta_+)\\
&\qquad\qquad-e^\kappa\underbrace{\int_{-\infty}^{\infty}\nu'(s)ds}_{=-1}\int_{0}^{1}\dot{\chi}(t)\max\Big\{\max_{x\in\Sigma}\Big(h_{\chi(t)}(x)-k_{\chi(t)}(x)\Big),0\Big\}dt\\
&=\mathcal{A}_{\varphi}(u_-,\eta_-)-\mathcal{A}_{\psi}(u_+,\eta_+)+e^\kappa\int_{0}^{1}\max\Big\{\max_{x\in\Sigma}\Big(h_t(x)-k_t(x)\Big),0\Big\}dt\\
&\leq\mathcal{A}_{\varphi}(u_-,\eta_-)-\mathcal{A}_{\psi}(u_+,\eta_+)+e^\kappa \max\big\{\| h-k\|_+,0\big\}.
\end{align}
This proves estimate \eqref{eq:the constant K}.

\section{\label{sec:Spectral-invariants}Spectral invariants and orderability}

Throughout this section we require Assumption \textbf{(A)} from the
Introduction to hold. More precisely, recall we say that a closed
connected coorientable contact manifold $(\Sigma,\xi)$ satisfies
Assumption \textbf{(A)} if:
\begin{description}
\item [Assumption {(A)}] $(\Sigma,\xi)$ admits a Liouville filling $(W,d\lambda)$
such that $\alpha:=\lambda|_{\Sigma}$ is Morse-Bott and the Rabinowitz
Floer homology $\RFH_{*}(\Sigma,W)$ is non-zero. \end{description}
\begin{Def}
\label{def of c}
Fix a non-zero class $ Z \in \RFH_*( \Sigma,W)$ and let $\varphi$ denote a non-degenerate path. We define
its \emph{spectral number} by
\begin{equation}
c(\varphi, Z):=\inf\left\{ a\in\mathbbm{R} \mid Z_{\varphi} \in j_{\varphi}^{a}(\RFH_{*}^{a}(\mathcal{A}_{\varphi},W))\right\},\label{eq:def of c}
\end{equation}
where we use the notation $ Z_\varphi$ from \eqref{eq:Z_varphi}.
\end{Def}
 Throughout the rest of the paper, the letter $ Z$ denotes a non-zero class in $ \RFH_*( \Sigma,W)$.
\begin{Prop}
\label{prop:the_comparison_prop}
Let $\varphi$ and $\psi$ be two non-degenerate paths. Then we have the estimate
\begin{align}
\label{eqn:estimate_1}
c(\psi, Z)&\leq c(\varphi , Z )+K(\varphi,\psi) \\
&\leq c(\varphi,Z)+ e^{\max\{\kappa(\varphi),\kappa(\psi)\}} \max\big\{\| h - k \|_+,0\big\},
\end{align}
where $h$ and $k$ are the contact Hamiltonians of $\varphi$ and $\psi$, respectively. In particular, we have
\begin{equation}
\label{eqn:estimate_2}
h_t(x)\leq k_t(x)\; \forall x\in\Sigma,t\in[0,1]\quad\Longrightarrow\quad c(\varphi,Z)\geq c(\psi,Z).
\end{equation}
\end{Prop}

\begin{proof}
This follows immediately from the definition of the spectral number together with \eqref{eq:fil cont hom} and the estimate \eqref{eq:the constant K}. 
\end{proof}

\begin{Lemma}
\label{lem: properties of c}
For any non-degenerate path $\varphi\in\mathcal{P}\Cont_{0}(\Sigma,\xi)$ the spectral numbers are all critical values of $\mathcal{A}_\varphi$, i.e.~$c(\varphi,Z)\in \mathrm{Spec}(\varphi)$. 

Moreover $c( \cdot, Z)$ admits a unique extension to all of $\mathcal{P}\Cont_{0}(\Sigma,\xi)$: given
a degenerate path $\varphi$, set 
\begin{equation}
c(\varphi,Z):=\lim_{k}c(\varphi_{k},Z),\label{eq:how defined in degen case-1}
\end{equation}
where $\varphi_{k}\rightarrow\varphi$ is any sequence of non-degenerate
paths converging to $\varphi$ in $C^2$. The extension still satisfies $c(\varphi,Z)\in \mathrm{Spec}(\varphi) $ and the estimates \eqref{eqn:estimate_1} and \eqref{eqn:estimate_2}. In particular, $c( \cdot, Z):\mathcal{P}\mathrm{Cont}_0(\Sigma,\xi) \rightarrow \mathbbm{R}$ is a continuous function when we equip $\mathcal{P}\mathrm{Cont}_0(\Sigma,\xi)$ with the $C^2$-topology. 
\end{Lemma}

\begin{proof}
The assertion $c(\varphi,Z)\in\mathrm{\mathrm{Spec}}(\varphi)$ follows immediately from the fact that $\RFH_{*}^{c}( \mathcal{A}_\varphi ,W)$ only changes if $c( \cdot, Z)$ crosses a critical value of $\mathcal{A}_\varphi$, compare the discussion below \eqref{proj_map}.

To prove the existence of the extension we are required to prove that the limit exists and is independent of the choice of approximating sequence $\varphi_k$. We denote by $h_k$ the corresponding contact Hamiltonians. Since we assume that $\varphi_k$ converges to $\varphi$ in $C^2$ it follows that $\kappa(\varphi_k)\to\kappa(\varphi)$ and $h_k\to h$, the contact Hamiltonian of $\varphi$. From Proposition \ref{prop:the_comparison_prop} we conclude that $(c(\varphi_k,Z))$ converges and in the same way independence of the approximating sequence $(\varphi_k)$ is proved. That  $c(\varphi,Z)\in\mathrm{Spec}(\varphi)$ and the estimates \eqref{eqn:estimate_1} and \eqref{eqn:estimate_2} hold follows from the definition of $c( \cdot, Z)$ as a limit.
\end{proof}

\begin{Lemma}
\label{lem:descends}
The map $c( \cdot, Z):\mathcal{P}\mathrm{Cont}_0(\Sigma,\xi) \rightarrow \mathbbm{R}$ descends to give a well defined map $c( \cdot, Z):\widetilde{\mathrm{Cont}}_0(\Sigma,\xi) \rightarrow \mathbbm{R}$.
\end{Lemma}

\begin{proof}
We recall from Remark \ref{Rmk: spec only depends on the end} that $\mathrm{Spec}(\varphi)\subset \R$ is nowhere dense and actually only depends on the endpoint $\varphi_1$ of the path $\varphi$. Moreover, Lemma \ref{lem: properties of c} implies that $c( \cdot, Z)$ is a continuous map. If we vary the path $\varphi$ while keeping the endpoints fixed the continuous map $c( \cdot, Z)$ takes values in the fixed, nowhere dense set $\mathrm{Spec}(\varphi_1)$, thus is constant. This proves the Lemma.
\end{proof}

\begin{Lemma}
\label{lem:c_reeb}
For any $ T \in \R$ one has 
\begin{equation}
  c( \theta^T ,Z ) = -T + c( \mathrm{id}_{ \Sigma}, Z).
\end{equation}
\end{Lemma}
\begin{proof}
One has $\mathrm{Spec}( \theta^T ) =  -T + \mathrm{Spec}(\mathrm{id}_{ \Sigma})$. Since $\mathrm{Spec}( \mathrm{id}_{ \Sigma})$ is nowhere dense and $c( \cdot, Z)$ is continuous the result follows.
\end{proof}

\begin{Rmk}
Proposition \ref{prop:the_comparison_prop} and Lemmata  \ref{lem: properties of c}, \ref{lem:descends},  \ref{lem:c_reeb} constitute Theorem \ref{Theorem-A} from the introduction.
\end{Rmk}

Given a path $\varphi$ of contactomorphisms, we define
the \emph{support} of $\varphi$, 
\begin{equation}
\mathfrak{S}(\varphi):=\bigcup_{0\leq t\leq1}\mathrm{supp}(\varphi_{t}),
\end{equation}
where $\mathrm{supp}(\varphi_{t}):=\overline{\left\{ x\in\Sigma\mid\varphi_{t}(x)\ne x\right\} }$. 

\begin{Def}
\label{def of cU}For an open set $U\subset\Sigma$ we set
\begin{equation}
c(U,Z):=\sup\left\{ c(\varphi,Z)\mid\varphi\in\widetilde{\Cont_{0}}(\Sigma,\xi),\ \mathfrak{S}(\varphi)\subset U\right\} \in(-\infty,\infty].
\end{equation}
\end{Def}

\begin{Ex}
\label{By-Example-infinit}By Lemma \ref{lem:c_reeb} one has immediately
that $c(\Sigma,Z)=\infty$ for any non-zero class $ Z$.
\end{Ex}

Recall from \eqref{eq:singular_homology_canonical} that  the fundamental class $ [ \Sigma]$ defines a class in $\RFH^{(- \varepsilon, \varepsilon)}_n( \Sigma,W)$. We now strengthen Assumption \textbf{(A)} as follows:
\begin{description}\label{ass:a+}
\item [Assumption {(A)$ ^+$}] $(\Sigma,\xi)$ admits a Liouville filling $(W,d\lambda)$
such that $\alpha:=\lambda|_{\Sigma}$ is Morse-Bott. Moreover there exists a non-zero class  $\mu_\Sigma \in \RFH_n( \Sigma,W)$ such that $  p^{ - \varepsilon}(  \mu_{ \Sigma} ) = j^{ \varepsilon}( [ \Sigma])$. See statement (2) \vpageref{enu:The-Rabinowitz-Floer2} for the definition of the maps $p^{-\varepsilon}$ and $j^{\varepsilon}$.
\end{description}
 
\begin{Rmk}\label{rmk:when_mu_Sigma_exisits}
Here are two instances where this assumption is satisfied. Firstly, if $  \Sigma$ is a unit cotangent bundle then such a non-zero class  $ \mu_{ \Sigma}$ exists; this follows from eg.~work of Abbondandolo-Schwarz \cite{AbbondandoloSchwarz2009}. Secondly, if $ ( \Sigma, \xi)$ admits a Liouville filling $(W, d \lambda)$ such that there are no contractible (in $W$) Reeb orbits, then again such a class $ \mu_{ \Sigma}$ trivially exists since in this case we can take $ \varepsilon = + \infty$ in \eqref{eq:singular_homology_canonical}, ie.~$\RFH_*( \Sigma,W) = \mathrm{H}_{*+n-1}( \Sigma; \Z_2)$.  
\end{Rmk}

\begin{Lemma}\label{lem:c(id,mu)=0}
It holds
\begin{equation}
c(\id_\Sigma,\mu_\Sigma)=0.
\end{equation}
\end{Lemma}

\begin{proof}
This follows from the definition of spectral numbers and the fact that $p^{ - \varepsilon}(  \mu_{ \Sigma} ) = j^{ \varepsilon}( [ \Sigma])$
\end{proof}

Let us show that $c(U)>0$ for any non-empty set $U\subset\Sigma$.
\begin{Prop} \label{prop:c(U)_is_positive}
Given any non-empty open set $U\subset\Sigma$,
there exists $\psi\in\widetilde{\mathrm{Cont}}_{0}(\Sigma,\xi)$
such that $\mathfrak{S}(\psi)\subset U$ and $c(\psi,\mu_\Sigma)>0$. 
\end{Prop}

\begin{proof}
We prove the proposition in three steps.

\textbf{Step 1. }We use an idea from Sandon \cite{Sandon2013}.
Fix a $C^{2}$-small function $b:\Sigma\rightarrow\mathbbm{R}$. We
use $b$ to build a contactomorphism $\Psi:T^{*}\Sigma\times\mathbbm{R}\rightarrow T^{*}\Sigma\times\mathbbm{R}$,
where the 1-jet bundle $T^{*}\Sigma\times\mathbbm{R}$ is equipped
with the standard contact form $\lambda_{0}+d\tau$ and $\lambda_{0}=pdx$
in local coordinates. Namely, we set 
\begin{equation}
\Psi(x,p,\tau)=(x,p-db(x),\tau+b(x)).
\end{equation}
Note that critical points of $b$ are in 1-1 correspondence with Reeb
chords between the two Legendrians $\Sigma\times\{0\}$ and $\Psi(\Sigma\times\{0\})$
(where $\Sigma\subset T^{*}\Sigma$ is the zero section). Since $b$
is assumed to be $C^{2}$-small, $\Psi$ determines a contactomorphism
of $\psi$ of $(\Sigma,\alpha)$, defined as follows. Firstly, Weinstein's
neighborhood theorem for Legendrian submanifolds (see \cite[Theorem 2.2.4]{AbbasHofer})
implies that there is an exact contactomorphism 
\begin{equation}
\Xi:N\times(-\delta,\delta)\rightarrow Q\times(-\varepsilon,\varepsilon)\label{eq:the contacto big xi-1}
\end{equation}
between a neighborhood $N\times(-\delta,\delta)$ of $\Sigma\times\{0\}$
inside $T^{*}\Sigma\times\mathbbm{R}$ and a neighborhood $Q\times(-\varepsilon,\varepsilon)$
of $\Delta\times\{0\}$ inside $\Sigma\times\Sigma\times\mathbbm{R}$,
where $\Delta$ is the diagonal in $\Sigma\times\Sigma$. Here $\Sigma\times\Sigma\times\mathbbm{R}$
is equipped with the contact form $e^{r}\mbox{pr}_{1}^{*}\alpha-\mbox{pr}_{2}^{*}\alpha$,
where $\mbox{pr}_{j}$ is the projection onto the $j$th factor. The
contactomorphism $\psi$ is then defined by looking at the restriction
of $\Xi\circ\Psi\circ\Xi^{-1}$ to $\Delta\times\{0\}$ inside $Q\times(-\varepsilon,\varepsilon)$;
we can write

\begin{equation}
\Xi\circ\Psi\circ\Xi^{-1}(x,x,0)=:(x,\psi(x),\log\rho(x)),
\end{equation}
for $\psi:\Sigma\rightarrow\Sigma$ and $\rho:\Sigma\rightarrow(0,\infty)$.
One readily checks that $\psi^{*}\alpha=\rho\alpha$, and hence $\psi$
is a contactomorphism. 

Similarly, if we start with an isotopy $\{b_{t}\}_{0\leq t\leq1}$
with $b_{0}=0$ then we obtain a path $\psi=\{\psi_{t}\}_{0\leq t\leq1}$
of contactomorphisms with $\psi_{0}=\id_{\Sigma}$. In this
case one can check that the contact Hamiltonian of $\psi$ is $-b_{t}$:
\begin{equation}
\alpha\left(\frac{\partial}{\partial t}\psi_{t}\right)=-b_{t}\circ\psi_{t}.
\end{equation}
The key point now is that the translated points $x\in\Sigma$ of $\psi_{1}$
with time-shift $\eta\in(-\varepsilon,\varepsilon)$ are in 1-1 correspondence
with the critical points of $b_{1}$: if $x\in\mathrm{Crit}(b_{1})$
then 
\begin{equation}
\psi(x)=\theta^{-b_{1}(x)}(x).
\end{equation}
Thus for each $x\in\mathrm{Crit}(b_{1})$ there is a critical point
$(u_{x},b_{1}(x))\in\mathrm{Crit}(\mathcal{A}_{\psi})$, and any critical
point $(u,\eta)$ of $\mathcal{A}_{\psi}$ not of this form necessarily
satisfies $\left|\eta\right|>\varepsilon$.

\textbf{Step 2. }Suppose now that we start with a $C^{2}$-small Morse-Bott
function $b$ on $\Sigma$. Define $b_{t}:=tb$ for $t\in[0,1]$,
and let $\psi=\{\psi_{t}\}_{0\leq t\leq1}$ denote the corresponding
path of contactomorphisms. Then if $x\in\mathrm{Crit}(b)$ then the
critical point $(u_{x},b(x))$ belongs to a Morse-Bott component of
$\mathcal{A}_{\psi}$, and moreover we claim that 
\begin{equation}
\mu(u_{x},b(x))=1-n+\mbox{ind}_{b}(x),
\end{equation}
where $\mbox{ind}_{b}(x)$ denotes the maximal dimension of a subspace
on which the Hessian $\mbox{Hess}_{b}(x)$ of $b$ at $x$ is \emph{strictly}
negative definite.

To see this, we consider the Hamiltonian diffeomorphism $\Phi$ of
$T^{*}\Sigma\times\mathbbm{R}\times\mathbbm{R}$ obtained by lifting
$\Psi$, which as $\Psi$ is exact, is given simply by 
\begin{equation}
\Phi(q,p,\tau,\sigma)=\left(\Psi(q,p,\tau),\sigma\right).
\end{equation}
A translated point $x$ of $\psi_{1}$ gives rise to the following
path of Lagrangian subspaces:
\begin{equation}
L_{t}:=\left\{ (\hat{x},-t\mbox{Hess}_{b}(x)(\hat{x}),0,\hat{\sigma})\mid\hat{x}\in T_{x}\Sigma,\,\hat{\sigma}\in\mathbbm{R}\right\} \subset T_{(x,0,0,\sigma)}(T^{*}\Sigma\times\mathbbm{R}\times\mathbbm{R}),
\end{equation}
The desired index is then given by
\begin{equation}
\mu(u_{x},b(x))=1-n+\mu_{\textrm{RS}}(L_{0},L_{1}),
\end{equation}
which in this case is just $1-n+\mbox{ind}_{b}(x)$ as claimed; note
that the $1-n$ summand comes from the normalization used in the definition
of the Rabinowitz index (\ref{eq:grading}) above, and we are using
the grading convention from Remark \ref{Our-grading-convention}.

\textbf{Step 3.} We now prove the theorem. Suppose $U\subset\Sigma$
is open and non-empty. Choose a function $b:\Sigma\rightarrow[0,\infty)$ such that
$\mbox{supp}(b)\subset U$ and such that $b$ is Morse on the interior
of its support. Moreover we insist that $b$ has a unique maximum
$x_{0}\in\Sigma$, with $0<b(x_{\max})<\varepsilon$, where $\varepsilon$
is as in (\ref{eq:the contacto big xi-1}).  Let $\psi$ be as in Step 2. Since the contact Hamiltonian of $\psi$ is $-b$ we can estimate
\begin{equation}
K(\id_{\Sigma},\psi)\leq e^{\kappa(\psi)} \tfrac12 b(x_{\max}).
\end{equation}
From Proposition \ref{prop:the_comparison_prop} and Lemma \ref{lem:c(id,mu)=0} we obtain 
\begin{align}
c(\psi,\mu_\Sigma)&\leq \underbrace{c(\id_{\Sigma},\mu_\Sigma)}_{=0}+K(\id_{\Sigma},\psi)\\
&\leq e^{\kappa(\psi)} \tfrac12 b(x_{\max}).
\end{align}
We now assume in addition that $e^{\kappa(\psi)} \tfrac12 b(x_{\max})<\epsilon$, too. Since the contact Hamiltonian $-b$ of $\psi$ is non-positive,
we have from \eqref{eqn:estimate_2} that
\begin{equation}
0\leq c(\psi,\mu_\Sigma)< \epsilon.
\end{equation}
We recall from Step 1 that any critical point $(u,\eta)$ of $\mathcal{A}_{\psi}$ which is not of the form $(u_x,b(x))$
satisfies $\left|\eta\right|>\varepsilon$. Thus $c(\psi,\mu_\Sigma)$ is necessarily a critical value of $b$. Since $\mu_\Sigma$ has index $n$, and $x_{\max}$ is the only critical point of $b$
of index $2n-1$ (so that the corresponding critical point $(u_{x_{\max}},b(x_{\max}))$
has index $1-n+2n-1=n$), we see that 
\begin{equation}
c(\psi,\mu_\Sigma)=b(x_{\max})>0.
\end{equation}
The proof is complete. 
\end{proof}

The following Corollary is Theorem \ref{thm:weak_inequality} from the Introduction.

\begin{Cor}
\label{cor:posi_but_not_negi}
Suppose $\varphi\in\widetilde{\mathrm{Cont}}_{0}(\Sigma,\xi)$
has contact Hamiltonian $h_{t}$. Assume $h_t\leq0$ and there exists $x\in\Sigma$ such that $h_t(x)<0$ for all $t\in[0,1]$. Then $c(\varphi,\mu_\Sigma)>0$.
\end{Cor}

\begin{proof} 
There exists a function $b:\Sigma\to[0,\infty)$ satisfying all the properties from the proof of Proposition \ref{prop:c(U)_is_positive}
and  in addition that
\begin{equation}
-tb(x)\geq h_t(x)\quad\forall x\in\Sigma,t\in[0,1].
\end{equation}
Let $ \psi = \{ \psi_t \}_{ 0 \le t \le 1}$ denote the contact isotopy whose contact Hamiltonian is $-tb$. Then Proposition \ref{prop:the_comparison_prop} and Proposition \ref{prop:c(U)_is_positive} imply that
\begin{equation}
0<b(x_{\max})=c(\psi,\mu_\Sigma)\leq c(\varphi,\mu_\Sigma).
\end{equation}

\end{proof}

\begin{Rmk}
\label{rem:reversing_prev_cor}
One might wonder whether the analogue of Corollary \ref{cor:posi_but_not_negi} continues to hold if instead we assume that $h_t$ is non-negative and not identically zero. In the non-compact setting discussed in Section \ref{sec:Relations-to-symplectic} we will see that this is false. See Remark \ref{rem:cptly_sprtd_rf} and Appendix \ref{app:A-'compactly-supported} for more information.
\end{Rmk}

\section{\label{sec:The-periodic-case - capacities} Contact capacities }

Let us now assume that $(\Sigma,\xi)$ satisfies Assumption \textbf{(B)} from the Introduction:
\begin{description}
\item [Assumption {(B)}] $(\Sigma,\xi)$ admits a Liouville filling $(W,d\lambda)$
such that the Rabinowitz Floer homology $\RFH_{*}(\Sigma,W)$
is non-zero and such that $\alpha:=\lambda|_{\Sigma}$ is periodic. 
\end{description}

As before $Z$ denotes a non-zero class in $\RFH_{*}(\Sigma,W)$.

\begin{Def}
\label{Def of c bar when resonant1}We define for $\varphi\in\widetilde{\Cont}_{0}(\Sigma,\xi)$
an integer $\overline{c}(\varphi,Z)$ by 
\begin{equation}
\overline{c}(\varphi,Z):=\left\lceil c(\varphi,Z)\right\rceil .
\end{equation}
 
\end{Def}
The reason periodicity is helpful is this function $\overline{c}(\cdot,Z)$
is \emph{conjugation invariant}. We will prove this shortly in Proposition
\ref{lem the use of periodicity} below, but to begin with we present
the following lemma. Recall from Definition \ref{def: nondegen} that
we say $\varphi$ is \emph{non-resonant }if $\mathrm{Spec}({\varphi})\cap\mathbbm{Z}=\emptyset$.
\begin{Lemma}
\label{lem:resonant and nondegen}Suppose $\varphi$ is both resonant
and degenerate with $c(\varphi,Z)\in\mathbbm{Z}$. Then there exists
$\varphi^{\nu}\rightarrow\varphi$ such that $\varphi^{\nu}$ is resonant
and non-degenerate such that for all $\nu$ sufficiently large one
has $c(\varphi^{\nu},Z)=c(\varphi,Z)$. \end{Lemma}
\begin{proof}
Start with any sequence $(\varphi^{\nu})$ of non-degenerate paths
such that $\varphi^{\nu}\rightarrow\varphi$. Since $\varphi$ is
resonant, for each $\nu\in\mathbbm{N}$ there exists a translated point
$x^{\nu}\in\Sigma$ for $\varphi^{\nu}$ with time-shift $\eta^{\nu}$
such that $\eta^{\nu}\rightarrow0$ and such that for all $\nu$ sufficiently
large one has 
\begin{equation}
c(\varphi^{\nu},Z)=c(\varphi,Z)+\eta^{\nu}.
\end{equation}
The sequence $\theta^{-\eta^{\nu}}\circ\varphi^{\nu}$ still converges
to $\varphi$, and it is easy to check that $\theta^{-\eta^{\nu}}\circ\varphi^{\nu}$
is still non-degenerate, and for all $\nu$ sufficiently large one
has that 
\begin{equation}
c(\theta^{-\eta^{\nu}}\circ\varphi^{\nu},Z)=c(\varphi^{\nu},Z)-\eta^{\nu}=c(\varphi,Z)
\end{equation}
since
\begin{equation}
\mathrm{Spec}(\theta^T\circ\varphi)=T+\mathrm{Spec}(\varphi)
\end{equation}
and $c(\cdot,Z)$ is continuous.
\end{proof}

The following is Proposition \ref{prop:c_conjugation_invariant_introduction} from the Introduction.

\begin{Prop}
\label{lem the use of periodicity}The function $\overline{c}(\cdot,Z):\widetilde{\mathrm{Cont}}_{0}(\Sigma,\xi)\rightarrow\mathbbm{Z}$
is conjugation invariant: if $\psi\in\mathrm{Cont}_{0}(\Sigma,\xi)$
and $\varphi\in\widetilde{\mathrm{Cont}}_{0}(\Sigma,\xi)$ then
\begin{equation}
\bar{c}(\psi\varphi\psi^{-1},Z)=\bar{c}(\varphi,Z).
\end{equation}
\end{Prop}

\begin{proof}
Assume firstly that $\varphi$ is non-resonant, that is, $\mathrm{Spec}({\varphi})\cap\mathbbm{Z}=\emptyset$
(see Definition \ref{def: nondegen}). Fix $\psi\in\Cont_{0}(\Sigma,\xi)$ and let $\psi_{s}\in\Cont_{0}(\Sigma,\xi)$ be a path connecting
$\id_{\Sigma}$ to $\psi$. Then we consider the map 
\begin{equation}
s\mapsto c(\psi_{s}\varphi\psi_{s}^{-1},Z).
\end{equation}
Proposition \ref{prop:the_comparison_prop} implies that this map is continuous.
Lemma \ref{lem:integral_critical_points_are_conjugation_invariant}
implies that $\mathrm{Spec}({\psi_{s}\varphi\psi_{s}^{-1}})\cap\mathbbm{Z}=\emptyset$
for all $s\in[0,1]$, and hence $\left\lceil c(\psi\varphi\psi^{-1},Z)\right\rceil =\left\lceil c(\varphi,Z)\right\rceil $
as required.

There are now three cases to consider. Suppose that $\varphi$
is resonant but that $c(\varphi,Z)\notin\mathbbm{Z}$. Suppose $\psi\in\Cont_{0}(\Sigma,\xi)$.
Then for $\varphi'$ non-resonant and sufficiently close to $\varphi$
we have 
\begin{equation}
\overline{c}(\psi\varphi\psi^{-1},Z)=\overline{c}(\psi\varphi'\psi^{-1},Z)=\overline{c}(\varphi',Z)=\overline{c}(\varphi,Z),
\end{equation}
where the second equality used the step above. The next case is
when $\varphi$ is resonant and non-degenerate, with $c(\varphi,Z)\in\mathbbm{Z}$.
As before, given $\psi\in\Cont_{0}(\Sigma,\xi)$ we choose a
path $\psi_{s}$ connecting $\id_{\Sigma}$ to $\psi$. The
key point now is that for any $s_{0}\in[0,1]$, if $(u_{s_{0}},\eta_{s_{0}})$
is a critical point of $\mathcal{A}_{\psi_{s_{0}}\varphi\psi_{s_{0}}^{-1}}$
with $\eta_{s_{0}}\in\mathbbm{Z}$ then $(u_{s_{0}},\eta_{s_{0}})$
is automatically non-degenerate by the last statement of Lemma \ref{lem:integral_critical_points_are_conjugation_invariant}.
It follows that there exists $\varepsilon>0$ such that 
\begin{equation}
\mathrm{Spec}({\psi_{s}\varphi\varphi_{s}^{-1}})\cap[c(\varphi,Z)-\varepsilon,c(\varphi,Z)+\varepsilon]=\{c(\varphi,Z)\},
\end{equation}
and the result follows as above. The final case is when $\varphi$
is both resonant and degenerate and $c(\varphi,Z)\in\mathbbm{Z}$. In
this case we employ Lemma \ref{lem:resonant and nondegen} to find
a sequence $\varphi^{\nu}\rightarrow\varphi$ such that $\varphi^{\nu}$
is both resonant, non-degenerate, and such that for all large $\nu$
one has $c(\varphi^{\nu},Z)=c(\varphi,Z)$. The argument above then implies
that for any $\psi\in\Cont_{0}(\Sigma,\xi)$ and for all $\nu$
sufficiently large, $c(\psi\varphi^{\nu}\psi^{-1},Z)=c(\varphi^{\nu},Z)$
is an integer. Since $c(\psi\varphi^{\nu}\psi^{-1},Z)\rightarrow c(\psi\varphi\psi^{-1},Z)$
the result follows. 
\end{proof}

\begin{Cor} \label{c_bar_Reeb}
One has $\bar{c}(t\mapsto\theta^{tT},Z)=\left\lceil -T+c(\id_\Sigma,Z)\right\rceil $
for any $T\in\mathbbm{R}$. \end{Cor}
\begin{proof}
Lemma \ref{lem:c_reeb}. 
\end{proof}
We now define $\overline{c}(U,Z)$ in the same way as $c(U,Z)$ was defined
in Definition \ref{def of cU}.
\begin{Def}
For an open set $U\subset\Sigma$ we define the \emph{contact capacity}
\begin{equation}
\overline{c}(U,Z):=\sup\left\{ \overline{c}(\varphi,Z)\mid\varphi\in\widetilde{\Cont}_{0}(\Sigma,\xi),\ \mathfrak{S}(\varphi)\subset U\right\} \in\mathbbm{Z}\cup\{\infty\}.
\end{equation}
\end{Def}

\begin{Rmk}
The notion of contact capacity was introduced by Sandon in \cite{Sandon2011b}. She was the first to discover a connection between translated points and orderability and other contact rigidity phenomena. 
\end{Rmk}

The following is Corollary \ref{prop:cu_co_are_invariant_under_contactos_and_monotone_wrt_inclusions_INTRODUCTION} from the Introduction.

\begin{Prop}
\label{prop:cu_co_are_invariant_under_contactos_and_monotone_wrt_inclusions}
For all $\psi\in\mathrm{Cont}_{0}(\Sigma,\xi)$, one has
\begin{equation}
\overline{c}(\psi(U),Z)=\overline{c}(U,Z).
\end{equation}
 \end{Prop}
\begin{proof}
Since 
\begin{equation}
\mathfrak{S}(\psi\varphi\psi^{-1})=\psi(\mathfrak{S}(\varphi)),
\end{equation}
we conclude from Proposition \ref{lem the use of periodicity} that
\begin{align}
\overline{c}(U,Z) & =\sup\left\{ \overline{c}(\varphi,Z)\mid\varphi\in\widetilde{\Cont}_{0}(\Sigma,\xi),\ \mathfrak{S}(\varphi)\subset U\right\} \nonumber \\
 & =\sup\left\{ \overline{c}(\varphi,Z)\mid\varphi\in\widetilde{\Cont}_{0}(\Sigma,\xi),\ \psi(\mathfrak{S}(\varphi))\subset\psi(U)\right\} \nonumber \\
 & =\sup\left\{ \overline{c}(\varphi,Z)\mid\varphi\in\widetilde{\Cont}_{0}(\Sigma,\xi),\ \mathfrak{S}(\psi\varphi\psi^{-1})\subset\psi(U)\right\} \nonumber \\
 & =\sup\left\{ \overline{c}(\psi\varphi\psi^{-1},Z)\mid\varphi\in\widetilde{\Cont}_{0}(\Sigma,\xi),\ \mathfrak{S}(\psi\varphi\psi^{-1})\subset\psi(U)\right\} \nonumber \\
 & =\sup\left\{ \overline{c}(\mu,Z)\mid\mu\in\widetilde{\Cont}_{0}(\Sigma,\xi),\ \mathfrak{S}(\mu)\subset\psi(U)\right\} \nonumber \\
 & =\overline{c}(\psi(U),Z).
\end{align}
\end{proof}

For completeness we recall Theorem \ref{cor:non_squeezing_INTRODUCTION} which is proved in the Introduction.

\begin{Thm}
\label{cor:non_squeezing}
Let $U\subset V\subset\Sigma$ be open sets and
assume that there exists $\varphi\in\mathrm{Cont}_{0}(\Sigma,\xi)$
with $\varphi(V)\subset U$. Then 
\begin{equation}
\overline{c}(U,Z)=\overline{c}(V,Z).
\end{equation}
In particular, if $\overline{c}(U,Z)<\overline{c}(V,Z)$ then there exists
no contact isotopy mapping $V$ into $U$. 
\end{Thm}

\begin{Rmk}
\label{rem:when_c_is_never_zero_on_open_sets}
If we assume that $(\Sigma,\xi)$ satisfies both assumption \textbf{(A)$^+$} and \textbf{(B)} then we know $\overline{c}(U,\mu_\Sigma)>0$ whenever $U\subset\Sigma$ is a nonempty open subset.  Unfortunately in general we do not know how to prove that $\bar{c}(U, Z)$ is ever finite. Nevertheless, in certain situations it \emph{is} possible to prove finiteness of the capacities, for instance when the subset $U$ is displaceable. In particular this is the case in the setting described in the next section, see Corollary \ref{cor:its_finite}. 
\end{Rmk}

\section{Prequantization spaces}
\label{sec:Relations-to-symplectic}

\subsection{\label{sub:Hamiltonian-Floer-homology}Hamiltonian Floer homology }

Fix a Liouville domain $(M_1,d\gamma_1)$. Let $S:=\partial M_1$
and $\kappa:=\gamma_1|_{S}$, so that $(S,\kappa)$ is a contact
manifold. Let $(M,d\gamma)$ denote the completion of $M_1$, so
that $M=M_1\cup_{S}(S\times[1,\infty))$. It is convenient in this
section to introduce the notation
\begin{equation}
M_\sigma:=\begin{cases}
M_1\backslash\big(S\times(\sigma,1)\big) & \text{if }0<\sigma<1,\\[.5ex]
M_1\cup_{S}\big(S\times[1,\sigma]\big) & \text{if }\sigma\geq1.
\end{cases}
\end{equation}
Note here we are using $\sigma$ to
denote the $\mathbbm{R}$-coordinate on the end of $M$ - this is so
as to avoid confusion in Section \ref{sub:Rabinowitz-Floer-homology on Mxs1},
when a second Liouville domain will come into play. 

Denote by $\mbox{Ham}_{c}(M,d\gamma)$ the group of Hamiltonian diffeomorphisms
$f$ on $M$ with compact support. As before, a \emph{path }$f=\{f_{t}\}_{0\leq t\leq1}$
of compactly supported Hamiltonian diffeomorphisms is assumed to be
smoothly parametrized and begin at the identity: $f_{0}=\id_{M}$.
Given such a path $f=\{f_{t}\}_{0\leq t\leq1}$, let $X_{f}$ denote
the time-dependent vector field on $M$ defined by 
\begin{equation}
\frac{\partial}{\partial t}f_{t}=X_{f_{t}}\circ f_{t}.
\end{equation}
The equation
\begin{equation}
f_{t}^{*}\gamma-\gamma=da_{t},\ \ \ a_{0}\equiv0\label{at-1}
\end{equation}
determines a smooth compactly supported function $a_{t}:M\to\mathbbm{R}$.
If we define
\begin{equation}
F_{t}=i_{X_{f_{t}}}\gamma-\left(\frac{\partial}{\partial t}a_{t}\right)\circ f_{t}^{-1},\label{eq:Ft from at-1}
\end{equation}
then $F_{t}$ generates $f_{t}$: $f_{t}=f_{F}^{t}$. We can recover
$a_{t}$ from $F_{t}$ via 
\begin{equation}
a_{t}=\int_{0}^{t}\left(i_{X_{f_{t}}}\gamma-F_{s}\right)\circ f_{s}ds\label{eq:relate a and F-1}
\end{equation}
(see for instance \cite[p294]{McDuffSalamon1998}). 

We briefly explain the construction of the Hamiltonian Floer homology
of $f$ in this section. The setting we consider here is a special
case of the one considered by Frauenfelder and Schlenk in \cite{FrauenfelderSchlenk2007},
to which we refer to for more details. However it will be convenient
for us to use the Morse-Bott framework developed by Frauenfelder \cite{Frauenfelder2004},
in order to make the link with the Rabinowitz Floer homology of $\Sigma:=M\times S^{1}$
clearer in the next section.

Let us first note that for a given $F\in C_{c}^{\infty}(S^{1}\times M,\mathbbm{R})$, the flow $f_{F}^{t}$
has many 1-periodic orbits, since $f_{F}^{t}$ is compactly supported. Of course, constant 1-periodic orbits outside the support of $f$ are uninteresting, and hence we introduce the following notation. Denote by 
\begin{equation}
\sigma(F):=\inf\left\{ \sigma>0\mid\mathfrak{S}(f_{F})\subseteq M_\sigma \right\} .
\end{equation}
Given a path $f=\{f_{t}\}_{0\leq t\leq1}$ in $\mbox{Ham}_{c}(M,d\gamma)$,
we set $\sigma(f):= \sigma(F)$, where $F$ is given
by (\ref{eq:Ft from at-1}). Next, we set 
\begin{equation}
\mathcal{P}_{F}:=\left\{ y\in M_{ \sigma(f)} \mid f_{F}^{1}(y)=y\right\} .
\end{equation}

\begin{Def}
Define a subset $\mathcal{H}_{c}^{\textrm{mb}}\subseteq C_{c}^{\infty}(S^{1}\times M,\mathbbm{R})$
(here the ``mb" stands for Morse-Bott) to consist of those functions
$F$ with the property that $\mathcal{P}_{F}$ is either a closed
submanifold of $M$ or an open domain whose closure is a compact manifold,
and for which
\begin{equation}
T_{y}\mathcal{P}_{F}=\ker(Df_{F}^{1}(y)-\mathbbm{1})\ \ \ \mbox{for all }y\in\mathcal{P}_{F}.\label{eq:morse bott}
\end{equation}
It is well known that the subset $\mathcal{H}_{c}^{\textrm{mb}}$
is generic in $C_{c}^{\infty}(S^{1}\times M,\mathbbm{R})$. We say
that a path $f=\{f_{t}\}_{0\leq t\leq1}$ is \emph{non-degenerate
}if the function $F$ defined in (\ref{eq:Ft from at-1}) belongs
to $\mathcal{H}_{c}^{\textrm{mb}}$. 
\end{Def}
We denote by $R_ \kappa$ the Reeb vector field of $\kappa$. Denote by
$\widehat{\mathcal{H}}$ the set of time-dependent smooth functions
$\widehat{F}$ on $M$ with the property that there exists $C>0$
such that $\widehat{F}_{t}|_{S\times[C,\infty)}$ is of the form $\widehat{F}_{t}(y,\sigma)=e(\sigma)$
for some smooth function $e:[C,\infty)\rightarrow\mathbbm{R}$ satisfying
\begin{equation}
\label{eq:asymptotic m behaviour-2}
0\leq e'(\sigma)<\wp(S,\kappa).
\end{equation}
Here 
\begin{equation}
\wp(S,\kappa):=\inf\left\{ T>0\mid\exists\ \text{a closed Reeb orbit of }R_ \kappa \text{ of period }T>0\right\} .
\end{equation}
This ensures that if $ $$\varphi_{\widehat{F}}^{1}$ denotes the
flow of $\widehat{F}$ then $\varphi_{\widehat{F}}^{1}$ has no non-constant
1-periodic orbits on $S\times(C,\infty)$. Note that if $F\in C_{c}^{\infty}(S^{1}\times M,\mathbbm{R})$
then one can find $\widehat{F}\in\widehat{\mathcal{H}}$ such that
$\widehat{F}|_{S^{1}\times M_{ \sigma(F)}}=F$. 

As a special case of this construction, consider a function $O$ on $M$ defined by requiring that $O=0$
on the interior $M_1^{ \circ}$ of $M_1$ and that 
\begin{equation}
O(y,\sigma)=e(\sigma)\label{eq:def of G}
\end{equation}
on $S\times [1, \infty)$, where $e(1)=0$ and $e$ satisfies \eqref{eq:asymptotic m behaviour-2}. 
In this case one has  
\begin{equation}
\label{eq:int M0}
\mathcal{P}_{O}=M_1,
\end{equation}
where points in $M_1$ are thought of as constant loops. In particular,
$f_{O}^{t}|_{M_1^ \circ }=\mathrm{id}_{M_1^\circ}$. Thus $ O$ is an extension of the zero function (generating the Hamiltonian diffeomorphism $\mathrm{id}_{ M_1}$ to $\widehat{ \mathcal{H}}$). 

\begin{Def}
Fix a path $f=\{f_{t}\}_{0\leq t\leq1}$, and let $F$ denote the
function defined in (\ref{eq:Ft from at-1}), and fix an extension
$\widehat{F}\in\widehat{\mathcal{H}}$ such that $\widehat{F}|_{S^{1}\times M_{ \sigma(F)}}=F$.
Recall that $\Lambda(M):=C_{\textrm{contr}}^{\infty}(S^{1},M)$. Define
the \emph{Hamiltonian action functional }$\mathcal{A}_{f}:\Lambda(M)\rightarrow\mathbbm{R}$
by 
\begin{equation}
\mathcal{A}_{f}(v):=\int_{0}^{1}v^{*}\gamma-\widehat{F}_{t}(v)dt.\label{A f}
\end{equation}

\end{Def}
Denote by $\mathrm{Crit}^{\circ}(\mathcal{A}_{f})$ the set of critical
points $v$ of $\mathcal{A}_{f}$ with $v(S^{1})\subseteq M_{ \sigma(f)}$.
Then $\mathrm{Crit}^{\circ}(\mathcal{A}_{f})$ doesn't depend on the
extension $\widehat{F}$ - in fact 
\begin{equation}
\mathrm{Crit}^{\circ}(\mathcal{A}_{f})=\mathcal{P}_{F},
\end{equation}
and hence the assumption (\ref{eq:morse bott}) implies that each
component of $\mathrm{Crit}^{\circ}(\mathcal{A}_{f})$ is a Morse-Bott
component for $\mathcal{A}_{f}$. 

Fix  $J\in\mathcal{J}_{\textrm{conv}}(M)$ (cf. (\ref{eq:convex type})).
We define an $L^{2}$-inner product $\left\langle \left\langle \cdot,\cdot\right\rangle \right\rangle _{J}$
on $\Lambda(M)$ as before (cf. (\ref{eq:l2 norm}), only this time
there is no $bb'$ term). We denote by $\nabla_{J}\mathcal{A}_{F}$
the gradient of $\mathcal{A}_{F}$ with respect to $\left\langle \left\langle \cdot,\cdot\right\rangle \right\rangle _{J}$.\emph{
}Pick a Morse function $g:\mathrm{Crit}^{\circ}(\mathcal{A}_{f})\rightarrow\mathbbm{R}$
and a Riemannian metric $\varrho$ on $\mathrm{Crit}^{\circ}(\mathcal{A}_{f})$
such that $(g,\varrho)$ is a Morse-Smale pair. In the case where
$\mathcal{P}_{F}$ is an open domain in $M$ whose boundary is a compact
manifold, $g$ must be chosen so that $\left\langle dg,\mathbf{n}\right\rangle <0$
on the boundary, where $\mathbf{n}$ is an outward pointing normal.
As before we define moduli spaces $\mathcal{M}_{v^{-},v^{+}}(\mathcal{A}_{f},g,J,\varrho)$
of gradient flow lines with cascades for critical points $v^{\pm}\in\mathrm{Crit}(g)$.
This time we grade $v\in\mathrm{Crit}(g)$ simply by $\mu(v):=\mu_{\textrm{CZ}}(v)+\mathrm{ind}_{g}(v)$,
where $\mu_{\mathrm{CZ}}(v)$ is the Conley-Zehnder index. A standard
convexity argument gives the necessary compactness needed to define
Floer homology - see Frauenfelder-Schlenk \cite{FrauenfelderSchlenk2007}. 

Given $-\infty<a<b<\infty$ denote by $\mbox{CF}_{*}^{(a,b)}(\mathcal{A}_{f},g):=\mathrm{Crit}_{*}^{(a,b)}(g)\otimes\mathbbm{Z}_{2}$,
where $\mathrm{Crit}_{*}^{(a,b)}(g)$ denotes the set of critical points
$v$ of $g$ with $a<\mathcal{A}_{f}(v)<b$. As before one defines
a boundary operator $\partial$ on $\mbox{CF}_{*}^{(a,b)}(\mathcal{A}_{f},g)$
by counting the elements of the zero-dimensional parts of the moduli
spaces $\mathcal{M}_{v^{-},v^{+}}(\mathcal{A}_{f},g,J,\varrho)$ for
$v^{-}\ne v^{+}$. We denote by $\mathrm{HF}_{*}^{(a,b)}(\mathcal{A}_{f})$
the associated homology, which as the notation suggests, is independent
of the auxiliary data $(g,J,\varrho)$. In fact, one can also show
it is also independent of the choice of path $f$. We abbreviate $\mathrm{HF}_{*}^{a}(\mathcal{A}_{f}):=\mathrm{HF}_{*}^{(-\infty,a)}(\mathcal{A}_{f})$
and $\mathrm{HF}_{*}(\mathcal{A}_{f}):=\mathrm{HF}_{*}^{(-\infty,\infty)}(\mathcal{A}_{f})$.
We denote the natural maps $\mathrm{HF}_{*}^{a}(\mathcal{A}_{f})\rightarrow\mathrm{HF}_{*}(\mathcal{A}_{f})$
by $j_{f}^{a}$ in the same way as before. Under our grading convention
explained in Remark \ref{Our-grading-convention}, there is a canonical
isomorphism 
\begin{equation}\label{eq:computng HF}
\mathrm{HF}_{*}(\mathcal{A}_{f})\cong\mathrm{H}_{n+*}(M_1,\partial M_1 ; \Z_2)\cong\mathrm{H}^{n-*}(M_1 ; \Z_2).
\end{equation}
See the proof of Lemma \ref{ex the unit for cM} below for one way
to see this. 

Next, the Floer homology $\mathrm{HF}_{*}(\mathcal{A}_{f})$ carries
the structure of a unital ring. The unit lives in degree $n$ according to our sign conventions, and under the isomorphism \eqref{eq:computng HF}, the unit corresponds to the fundamental class $ [M_1] \in \mathrm{H}_{2n}(M_1 , \partial M_1 ; \Z_2)$; see Lemma \ref{ex the unit for cM} below. We denote the unit by $\mathbf{1}_{f}\in\mathrm{HF}_{n}(\mathcal{A}_{f})$. Since $\mathrm{HF}_{*}(\mathcal{A}_{f})$ is necessarily
non-zero, as before we can define the \emph{spectral number }
\begin{equation}
c_{M}(f):=\inf\left\{ a\in\mathbbm{R}\mid\mathbf{1}_{f}\in j_{f}^{a}(\mathrm{HF}_{*}^{a}(\mathcal{A}_{f}))\right\} .
\end{equation}
As before, $c_{M}$ is a well defined function 
\begin{equation}
c_{M}:\widetilde{\mbox{Ham}}_{c}(M,d\gamma)\rightarrow\mathbbm{R}.
\end{equation}
We can use $c_{M}$ to define a capacity on open subsets $\mathcal{O}\subset M$,
\begin{equation}
c_{M}(\mathcal{O}):=\sup\left\{ c_{M}(f)\mid\mathfrak{S}(f)\subset\mathcal{O}\right\} ,
\end{equation}
in the same way as before. We use the subscript $c_{M}$ to differentiate
it from the function $c$ associated to $\Sigma:=M\times S^{1}$ that
we will define shortly. 

\begin{Lemma}
\label{ex the unit for cM}
In the case of $ \id = \id_{M_1}$ the unit $\mathbf{1}= \mathbf{1}_{\id}$ is simply given by the fundamental class $[M_1]$, and thus $c_{M}(\id)=0$. 
\end{Lemma}

\begin{proof}
We define $\mathcal{A}_{\id}$ using the function $O$
defined in \eqref{eq:def of G}. Thus $\mathrm{Crit}^{\circ}(\mathcal{A}_{\id})=M_1$,
and every element of $\mathrm{Crit}^{\circ}(\mathcal{A}_{\id})$
has action value zero. Thus there are no gradient flow lines of $\mathcal{A}_{\id}$,
and hence the Floer complex $\mathrm{CF}_{*}(\mathcal{A}_{\id},g)$
reduces to the Morse complex of a Morse function $g$ on $M_1$.
Such a Morse function $g$ can be chosen so that $g>1$ on $M_1^ \circ $
and such that $g$ is the restriction of a Morse function $\widehat{g}:M\rightarrow\mathbbm{R}$
such that $\widehat{g}(y,\sigma)=\frac{1}{\sigma}$ on $S\times[1,\infty)$.
Thus this shows that 
\begin{equation}
\label{it_morse_hom}
\mathrm{HF}_{*}(\mathcal{A}_{\id})\cong\mathrm{HM}_{n+*}(g)\cong\mathrm{H}_{n+*}(M_1,\partial M_1; \Z_2),
\end{equation}
which proves \eqref{eq:computng HF}. 

It is possible to prove directly using Morse-Bott techniques that the isomorphisms in \eqref{it_morse_hom} are ring maps, and thus the  unit in $\mathrm{HF}_*( \mathcal{A}_{\id })$ is exactly the unit in Morse homology for $g$. The latter is of course the fundamental class $[M_1]$ under the isomorphism of the Morse homology of $g$ with
the relative homology of $(M_1,\partial M_1)$. In this situation however, we can simply make a degree argument: if the Morse function $g$ has a unique maximum at a point $y_{ \max}$ in $M_1^ \circ$ then one necessarily has that $[y_{\max}]$ is a cycle in $\mathrm{HF}_{n}(\mathcal{A}_{\id})$, and that fact
$\mathrm{HF}_{n}(\mathcal{A}_{\id})=\mathbbm{Z}_{2}[y_{\max}]$. Since the unit lives in degree $n$, it must therefore be precisely $[y_{ \max}]$.  
\end{proof}

\subsection{The prequantization space $\Sigma=M\times S^{1}$ }

The \emph{prequantization space } of
$M$ is the contact manifold $\Sigma:=M\times S^{1}$, equipped with
the contact structure $\xi:=\ker\,\alpha$, where
\begin{equation}
\alpha:=\gamma+d\tau,
\end{equation}
and $\tau$ is the coordinate on $S^{1}\cong\mathbbm{R}/\mathbbm{Z}$.
The last class of contact manifolds we study in this paper are these
prequantization spaces, which for convenience we refer to as Assumption
\textbf{(C)}:
\begin{description}
\item [Assumption {(C)}] $(\Sigma,\xi=\ker\,\alpha)$ is a prequantization space $\Sigma=M\times S^{1}$,
where $(M,d\gamma)$ is a Liouville manifold, and $\alpha=\gamma+d\tau$. 
\end{description}
In this case $\Sigma$ is obviously periodic, but it is \emph{not
}Liouville fillable in the previous sense. Aside from anything else,
$\Sigma$ is necessarily \emph{non-compact}. However $\Sigma$ does
still retain enough of the properties needed above in order to define
a Rabinowitz Floer homology, as will explain in the next section.

Let us denote by $\mathrm{Cont}_{0,c}(\Sigma,\xi)$ those contactomorphisms
$\varphi$ with compact support. There is a natural way to obtain
a path $\varphi=\{\varphi_{t}\}_{0\leq t\leq1}$ of compactly supported
contactomorphisms on $\Sigma$ from a path $f=\{f_{t}\}_{0\leq t\leq1}$
of compactly supported Hamiltonian diffeomorphisms on $M$. Indeed,
given such a path $f$, define $\varphi_{t}:\Sigma\rightarrow\Sigma$
by
\begin{equation}
\varphi_{t}(y,\tau):=\bigl(f_{t}(y),\underset{\textrm{mod }1}{\underbrace{\tau-a_{t}(y)}}\bigr),\label{eq:phi from f-1}
\end{equation}
where $a_{t}$ was defined in (\ref{at-1}). One easily checks that
$\varphi_{t}$ is an \emph{exact} contactomorphism. We say that the
contact isotopy $\varphi$ is the \emph{lift }of the Hamiltonian isotopy
$f$. In this case the contact Hamiltonian $h_{t}$ associated to
$\varphi_{t}$ is simply $F_{t}$: 
\begin{equation}
h_{t}\circ\varphi_{t}=\alpha\left(\frac{\partial}{\partial t}\varphi_{t}\right)=F_{t}\circ\varphi_{t},\label{eq:h equals F}
\end{equation}
where $F_{t}$ was defined in (\ref{eq:Ft from at-1}). Fix a function
$\widehat{F}\in\widehat{\mathcal{H}}$ such that $\widehat{F}=F$
on $S^{1}\times M_{ \sigma(F)}$, and define $\widehat{H}_{t}:S\Sigma\rightarrow\mathbbm{R}$
by $\widehat{H}_{t}:=r\widehat{F}_{t}$.

Consider again the Rabinowitz action functional $\mathcal{A}_{\varphi}:\Lambda(S\Sigma)\times\mathbbm{R}\rightarrow\mathbbm{R}$
defined as in (\ref{eq:Rab functional}), using $\widehat{H}_{t}$.
Suppose $(u,\eta)\in\mathrm{Crit}(\mathcal{A}_{\varphi})$. Write $u(t)=(v(t),\tau(t),r(t))\in M\times S^{1}\times(0,\infty)$.
Then from (\ref{eq:critical points}) we have 
\begin{align}
(f_{1}(v\left(\tfrac{1}{2}\right),\underset{\textrm{mod 1}}{\underbrace{\tau\left(\tfrac{1}{2}\right)-a_{1}\left(v\left(\tfrac{1}{2}\right)\right)}}\bigr) & =\varphi_{1}\left(u\left(\tfrac{1}{2}\right)\right),\nonumber \\
 & =\theta^{-\eta}\left(v\left(\tfrac{1}{2}\right),\tau\left(\tfrac{1}{2}\right)\right)\nonumber \\
 & =\bigl(v\left(\tfrac{1}{2}\right),\underset{\textrm{mod 1}}{\underbrace{\tau\left(\tfrac{1}{2}\right)-\eta}}\bigr)
\end{align}
and hence if $y:=v(\tfrac{1}{2})$ then $f_{1}(y)=y$ and $a_{1}(y)=\eta$
mod 1. Moreover since $\varphi_{t}$ is exact one has $r(t)\equiv1$
for all $t$ (cf. the last statement of Lemma \ref{Lemma:critical points}).
Since we only consider \emph{contractible }critical points of $\mathcal{A}_{\varphi}$,
we deduce: 
\begin{Lemma}
\label{lem:relating the two functionals}There exists a bijective
map 
\begin{equation}
\pi:\mathrm{Crit}(\mathcal{A}_{\varphi})\rightarrow\mathrm{Crit}(\mathcal{A}_{f})
\end{equation}
given by 
\begin{equation}
\pi(u=(v,\tau,r),\eta):=\left(t\mapsto f_{t}\left(v\left(\tfrac{1}{2}\right)\right)\right).
\end{equation}
Moreover 
\begin{equation}
\mathcal{A}_{\varphi}(u,\eta)=\mathcal{A}_{f}(\pi(u,\eta)).
\end{equation}
In particular, every critical point $(u,\eta)$ of $\mathcal{A}_{\varphi}$
has 
\begin{equation}
u(S^{1})\subseteq M_{ \sigma(f)} \times S^{1}\times\{1\}.
\end{equation}
\end{Lemma}

Given a contactomorphism $\varphi\in\mathrm{Cont}_{0,c}(\Sigma,\xi)$,
we denote by 
\begin{equation}
\sigma( \varphi) =\inf\left\{ \sigma>0\mid\mathfrak{S}(\varphi)\subseteq M_\sigma \times S^{1}\right\} .
\end{equation}
Thus if $\varphi$ is the lift of $f$ then 
\begin{equation}
\sigma( \varphi) = \sigma (f).
\end{equation}

\subsection{\label{sub:Rabinowitz-Floer-homology on Mxs1}Rabinowitz Floer homology
on $\Sigma$}

Let $P_1$ denote a torus with a disc removed, so that $\partial P_1=S^{1}$.
Equip $P_1$ with an exact symplectic form $d\beta_1$ such that
$\beta_1|_{\partial P_1}=d\tau$. Denote by $(P,d\beta)$ the
completion of $P$, so that 
\begin{equation}
\beta=rd\tau \quad \mathrm{on} \qquad \partial P_1\times[1,\infty).
\end{equation}
Consider 
\begin{equation}
W:=M\times P,
\end{equation}
equipped with the symplectic form $d\lambda$ where $\lambda:=\gamma+\beta$.
Note that $W$ is \emph{not} a Liouville filling of $\Sigma$. Indeed, firstly $W$ is not compact, and moreover when equipped with the symplectic form $d\lambda$, there
is no embedding $(S\Sigma,d(r\alpha))\hookrightarrow(W,d\lambda)$
that we can use in order to extend the Rabinowitz action functional
$\mathcal{A}_{\varphi}$ to a functional defined on all on $\Lambda(W)\times\mathbbm{R}$.
Ideally we would like to work with a symplectic form $\omega'$ such
that $\omega'|_{M\times\partial P_1 \times(0,\infty)}=d(r\alpha)$.
Sadly no such symplectic form exists, and even if one did it would
not have the ``right" properties at infinity. 

To circumvent this problem we will work with a family $ \{ \omega_s =  d\lambda_s \}_{s\geq1}$ of symplectic forms which satisfy: 

\begin{itemize}
	\item If $\varphi = \{ \varphi_t \}_{0 \le t \le 1}$ is any path of compactly support contactomorphisms on $\Sigma$ then there exists $s_0( \varphi) > 0$ with the property that for all $s > s_0( \varphi)$ the 1-form $\lambda_s $ agrees with $r \alpha$ on $\mathfrak{S}(\varphi)$.
	\item For every $s \ge 1$ the symplectic form $\omega_s$ is split-convex at infinity.
\end{itemize}

More precisely, we prove:
\begin{Lemma}
There exists a family $\{\lambda_s\}_{s\geq1} \subset \Omega^1(W)$ of 1-forms such that for all $s \ge 1$:

\begin{enumerate}
\item $\omega_s := d\lambda_s$ is a symplectic form on $W$.
\item Define
\begin{align}
W^+_s  & := \left(M\backslash M_{2s-1} \times P\right)\cup\left(M\times P\backslash P_{2s-1}\right), \label{W^+_s}\\
W^-_s & := M_s \times P_{\sfrac{1}{(2s-1)}}.
\end{align}
Then  
\begin{equation}
\lambda_s|_{W^+_s} = (2s-1)\gamma + \beta, \qquad \lambda_s|_{W^-_s} = \frac{1}{2s-1}\gamma + \beta.
\end{equation}
Thus $\omega_s$ is split-convex at infinity, and hence we can achieve compactness, see the proof of Theorem \ref{thm:rfh doesnt depend on s} below and also \cite{FrauenfelderSchlenk2007}. Moreover $\lambda_1 = \lambda$ everywhere.
\item For $s>1$, define 
\begin{equation}
V_s := M_s \times S^1 \times ( \sfrac{1}{s}, s) \subset S \Sigma.
\end{equation}
Then for each $s>1$, the natural embedding
\begin{equation}
\label{eq:embedding i}
\iota_{s}:V_s \hookrightarrow W
\end{equation}
 satisfies $\iota_s^*\lambda_s=r\alpha$.
\end{enumerate}
\end{Lemma}

\begin{proof}
Define a family $\{ f_s \}_{ s \ge1 }$ of smooth functions, see Figure \ref{fig_irida} :
\begin{equation}
f_s : [0,\infty)\times[0,\infty)\rightarrow(0,\infty)
\end{equation}
such that
\begin{equation}
f_s( \sigma, r) =
\begin{cases}
r, & ( \sigma, r) \in [0,s) \times (1/s,s), \\
\tfrac{1}{2s-1}, & ( \sigma, r) \in [0,s) \times (0,\tfrac{1}{2s-1}), \\
2s-1, & ( \sigma, r) \in [0,s) \times (2s-1, \infty), \\
2s-1, & ( \sigma, r) \in [s, \infty) \times (0,  \infty). \\
\end{cases}
\end{equation}
and finally such that 
\begin{equation}
\label{eq:why_symplectic}
\frac{\partial f_s}{\partial\sigma}(\sigma,r)\geq0, \qquad \text{for all } (s,\sigma,r)\in[1,\infty)\times[0,\infty)\times[0,\infty).
\end{equation}
The fact that such functions $f_s$ exist is clear from Figure \ref{fig_irida}.

\begin{figure}[ht]
\def\svgwidth{60ex} 
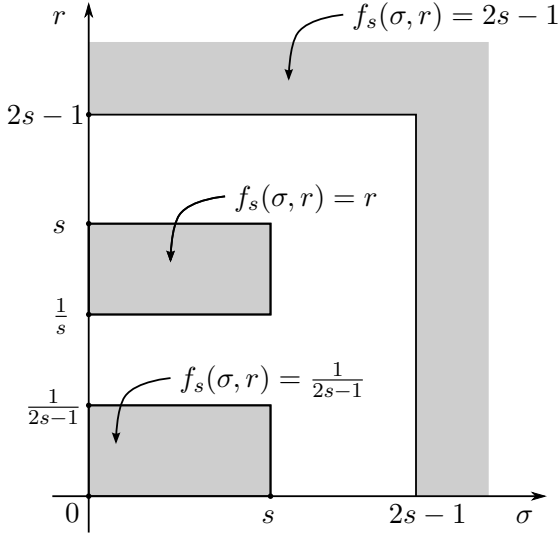
\caption{The function $f_s(\sigma, r)$}\label{fig_irida}
\end{figure}


On $M \setminus M_0 \times P \setminus P_0 $, where both the $\sigma$ and $r$-coordinates are defined, we set 
\begin{equation}
\lambda_s : =f_s \gamma + \beta.
\end{equation}
The condition \eqref{eq:why_symplectic} guarantees that $\omega_s := d \lambda_s$ is symplectic where defined, and it is clear that statements (2) and (3) from the Lemma are satisfied. It remains to extend $\lambda_s$ to all of $W$. This is done simply by ``continuity'': 
\begin{equation}
\lambda_s = \begin{cases}
f_s(0, r)\gamma + \beta, & \text{on } M_0 \times P \setminus P_0, \\
f_s(\sigma, 0)\gamma + \beta, & \text{on } M \setminus M_0 \times P_0, \\
\frac{1}{2s-1}\gamma + \beta, & \text{on } M_0 \times P_0.
\end{cases}
\end{equation}
\end{proof}

\begin{Def}
Given a path $\varphi = \{ \varphi_t \}_{0 \le t \le 1} $ of compactly supported contactomorphisms, we define the number $s_0 (\varphi) \ge 1$ by:
\begin{equation}
s_0 ( \varphi) := \inf \{ s \ge 1 \mid \mathfrak{S}( \varphi) \subset \iota_s(V_s) \},
\end{equation}
where $\iota_s$ is the embedding \eqref{eq:embedding i}.
\end{Def}

We now prove the following result. 
\begin{Thm}
\label{thm:rfh doesnt depend on s}
For any non-degenerate path $\varphi$
if $s >  s_{0}(\varphi)$
then it is possible to define the Rabinowitz Floer homology $\RFH_*(\mathcal{A}_{\varphi},W,\omega_s)$
(here the notation indicates that we are working with the symplectic
structure $\omega_s$ on $W$). Moreover the Rabinowitz Floer
homology is independent of the choice of $s >  s_0 (\varphi)$. 
\end{Thm}

\begin{proof}
Let $\mathcal{J}_{s}(W)$ denote the set of time-dependent almost
complex structures $J=\{J_{t}\}_{t\in S^{1}}$ on $W$ that are $\omega_s$
compatible, and satisfy: 
\begin{enumerate}
\item If $\iota_{s}$ is the embedding \eqref{eq:embedding i} then $\iota_{s}^{*}J\in\mathcal{J}_{\textrm{conv}}(V_{s}\subset S\Sigma)$. 
\item The restriction of $J$ to the subset $W^+_s$ defined in \eqref{W^+_s}
is \emph{split} - that is, there exist almost complex structures $J'\in\mathcal{J}_{\mathrm{conv}}(M)$
and $J''\in\mathcal{J}_{\mathrm{conv}}(P)$ such that $J=J' \oplus J''$
on this set. 
\end{enumerate}
Extend $\mathcal{A}_{\varphi}$ to a functional $\mathcal{A}_{\varphi}^{\kappa}$
defined on all of $\Lambda(W)\times\mathbbm{R}$ in the same way as
before, by replacing $\widehat{H}_{t}$ with a truncated function
$\widehat{H}_{t}^{\kappa}$ as in (\ref{eq:truncating H}). As with
the Hamiltonian Floer homology, we are now only interested in the
set $\mathrm{Crit}^{\circ}(\mathcal{A}_{\varphi}^{\kappa})$ of critical
points $(u=(v,\tau,r),\eta)$ of $\mathcal{A}_{\varphi}^{\kappa}$
with $v(S^{1})\subset M_{\sigma(\varphi)}$. For $s$ large
enough, all elements of $\mathrm{Crit}^{\circ}(\mathcal{A}_{\varphi}^{\kappa})$
are contained in $\iota(V_{s})$ and $\widehat{H}_{t}^{\kappa}$ is
constant outside on $W^+_{2s-1+\varepsilon}$ for some small $\varepsilon>0$.
Thus if we work with an almost complex structure $J\in\mathcal{J}_{s}(W)$,
the maximum principle prohibits the cylinder part of flow lines of $-\nabla_{J}\mathcal{A}_{\varphi}^{\kappa}$ from
ever entering $W^+_{2s-1+\varepsilon}$, see for instance \cite{FrauenfelderSchlenk2007}. Thus the Rabinowitz Floer
homology is well defined for this $s$. We point out that, since the cylinder part of flow lines stay in a compact subset of $W$, $L^\infty$-bounds on the Lagrange multiplier are derived as in \cite[Theorem 3.1]{CieliebakFrauenfelder2009}.

In order to prove independence of $s$, first note that for $s > \max\{s_{0}(\varphi),s_{0}(\psi)\}$
the continuation maps from points (1)-(4) \vpageref{enu:The-Rabinowitz-Floer1} 
show that 
\begin{equation}
\label{eq:indep of s 1}
\RFH_*(\mathcal{A}_{\varphi},W,\omega_s)\cong\RFH_*(\mathcal{A}_{\psi},W,\omega_s).
\end{equation}
Next we note that if $\id := \id_{M_1\times S^{1}}$ is the contactomorphism
with contact Hamiltonian $O$ as defined in \eqref{eq:def of G} then $s_0( \id)=1$. More generally, this is true for any \emph{exact }path $\varphi$ of contactomorphisms, since
in this case for any $\varepsilon>0$, every critical point of $\mathcal{A}_{\varphi}^{\varepsilon}$
is contained in $\Sigma\times\{0\}$ - see Lemma \ref{prop:linfinity-1}.
Thus by  \eqref{eq:indep of s 1} it suffices to show that $\RFH_*(\mathcal{A}_{\id},W,\omega_s)$
is independent of $s > 1$. But this is clear, since every critical point
of the Rabinowitz action functional $\mathcal{A}_{\id}$
has action value zero, as we are only looking at contractible critical
points and we have filled $S^{1}$ with a punctured torus $P_1$
rather than a disc $D^{2}$. Thus $\mathrm{Crit}^{\circ}(\mathcal{A}_{\id})=M_1\times S^1 \times \{ 0 \}$,
and hence regardless of which symplectic structure we use, as in Lemma
\ref{ex the unit for cM}, the Rabinowitz complex reduces to the Morse
complex of a Morse function $\tilde{g}:M_1\times S^{1}\rightarrow\mathbbm{R}$. In particular, it does not depend on $s$. 
\end{proof}

We denote by $\RFH_{*}(\mathcal{A}_{\varphi},W)$ the groups
$\RFH_{*}(\mathcal{A}_{\varphi},W,\omega_s)$ for any $s >  s_{0}(\varphi)$.

\begin{Thm}
\label{thm computation}
If $\varphi=\{\varphi_{t}\}_{0\leq t\leq1}$
is the lift of $f=\{f_{t}\}_{0\leq t\leq1}$ then there exists a natural
isomorphism 
\begin{equation}
\RFH_{*}(\mathcal{A}_{\varphi},W)\cong\HF_{*}(\mathcal{A}_{f})\otimes\mbox{\emph{H}}_{*}(S^{1}; \Z_2).
\end{equation}
\end{Thm}

\begin{proof}
By naturality it suffices to prove the theorem in the case $f=\id_{M_1}$ and
$\varphi= \id := \id_{M_1\times S^{1}}$. In this case as in the proof
of the last part of Theorem \ref{thm:rfh doesnt depend on s}, one
has
\begin{equation}
\RFH_{*}(\mathcal{A}_{\id},W)\cong\HM_{*+n}(\tilde{g}),
\end{equation}
where $\tilde{g}$ is a Morse function on $M_1 \times S^{1}$. We
choose $\tilde{g}=(g,g')$, where $g$ is the Morse function considered
in the proof of Lemma \ref{ex the unit for cM}, and $g':S^{1}\rightarrow\mathbbm{R}$
is a Morse function with two critical points $\tau_{\min}$ and $\tau_{\max}$.
This gives 
\begin{align}
\RFH_{*}(\mathcal{A}_{\id},W) & \cong\HM_{n+*}(\tilde{g}) \\ 
& \cong\HM_{*}(g)\otimes\HM_{*}(g') \\
& \cong \HM_{n+*}(M_1,\partial M_1)\otimes \mathrm{H}_*(S^{1}; \Z_2).
\end{align}

This completes the proof.
\end{proof}

\begin{Rmk}\label{rmk:mu_Sigma_prequant}
As pointed out in Remark \ref{rmk:when_mu_Sigma_exisits} there exists a certain non-zero class $\mu_\Sigma\in\RFH(\Sigma,W)$. The image of $\mu_\Sigma$ under the isomorphisms from Theorem \ref{thm computation} and \eqref{eq:computng HF} is the class $[M_1]\otimes[S^{1}]$.
\end{Rmk}

\subsection{\label{sub:Applications}Relating the capacities}

\begin{Def}
 We define $c(\varphi):=c(\varphi,\mu_\Sigma)$ in the same way as before for $\varphi\in\widetilde{\mathrm{Cont}}_{0,c}(\Sigma,\xi)$.
\end{Def}

As long as we work with compactly supported contactomorphisms Proposition \ref{prop:the_comparison_prop} remains true and its prove is literally the same.
 
\begin{Prop}\label{prop:the_comparison_prop_prequant}
Let $\varphi,\,\psi\in\widetilde{\mathrm{Cont}}_{0,c}(\Sigma,\xi)$ be two non-degenerate paths. Then we have the estimate
\begin{align}
\label{eqn:estimate_11}
c(\psi)&\leq c(\varphi)+K(\varphi,\psi) \\
&\leq c(\varphi)+ e^{\max\{\kappa(\varphi),\kappa(\psi)\}} \| h - k \|_+,
\end{align}
where $h$ and $k$ are the contact Hamiltonians of $\varphi$ and $\psi$, respectively. In particular, we have
In particular, we have
\begin{equation}
\label{eqn:estimate_21}
h_t(x)\leq k_t(x)\; \forall x\in\Sigma,t\in[0,1]\quad\Longrightarrow\quad c(\varphi)\geq c(\psi)
\end{equation}
and the same implication with nonstrict inequalities.
\end{Prop}

The analogue of Corollary \ref{cor:posi_but_not_negi} remains true, too, again with the same proof.

\begin{Cor}
\label{cor:posi_but_not_negi_prequant}
Suppose $\varphi\in\widetilde{\mathrm{Cont}}_{0,c}(\Sigma,\xi)$
has contact Hamiltonian $h_{t}$. Assume $h_t\leq0$ and $h_t\neq0$ for all $t\in[0,1]$. Then $c(\varphi)>0$.
\end{Cor}

\begin{Rmk}
\label{rem:cptly_sprtd_rf}
Recall in the closed case we proved that $c( t \mapsto \theta^{tT} ,Z) = -T + c( \id_{\Sigma}, Z)$ for any $T \in \R$ (cf. Statement (2) of Theorem \ref{Theorem-A}. In this setting the Reeb flow $\theta^t$ is of course \emph{not} compactly supported, and thus its spectral value is not defined. Nevertheless it is still possible to define a ``compactly supported Reeb flow'' $\vartheta^t : \Sigma \to \Sigma$ which agrees with the normal Reeb flow on a neighborhood of a given closed Reeb orbit.  For small $T$ it is still possible to compute the spectral numbers $c(\vartheta^T)$, but it is no longer the case that $c( \vartheta^T) = -T$. Indeed, whilst for negative $T$ one still has $c(\vartheta^T)= -T$, for positive $T$ one has $c(\vartheta^T) = 0$. This shows that Corollary \ref{cor:posi_but_not_negi_prequant} fails if one instead assume $h_t \ge 0$. Details are contained in Appendix \ref{app:A-'compactly-supported}.
\end{Rmk}

\begin{Def}
For an open non-empty set $U\subset\Sigma$ with compact closure we set
\begin{equation}
c(U):=\sup\left\{ c(\varphi)\mid\varphi\in\widetilde{\Cont_{0,c}}(\Sigma,\xi),\ \mathfrak{S}(\varphi)\subset U\right\} \in(-\infty,\infty].
\end{equation}
and
\begin{equation}
\overline{c}(U):=\lceil c(U)\rceil.
\end{equation}
\end{Def}

\begin{Thm}
\label{thm:equality of capacities-1}
Suppose $f\in\widetilde{\mbox{\emph{Ham}}}_{c}(M,d\gamma)$,
and let $\varphi\in\widetilde{\Cont}_{0,c}(\Sigma,\xi)$ denote
the lift of $f$. Then 
\begin{equation}
c_{M}(f)=c(\varphi).
\end{equation}
Moreover, if $\mathcal{O}\subset M$ is open with compact closure then 
\begin{equation}
c_{M}(\mathcal{O})=c(\mathcal{O}\times S^{1}).
\end{equation}
\end{Thm}

\begin{proof}
The first statement follows from Lemma \ref{lem:relating the two functionals}, Theorem \ref{thm computation}, and Remark \ref{rmk:mu_Sigma_prequant}. Thus clearly $c_{M}(\mathcal{O})\leq c(\mathcal{O}\times S^{1})$.
In order to complete the proof, we must show that given any $\psi\in\widetilde{\mathrm{Cont}}_{0,c}(\Sigma,\xi)$
with $\mathfrak{S}(\psi)\subset\mathcal{O}\times S^{1}$ there exists
$f\in\widetilde{\mbox{Ham}}_{c}(M,d\gamma)$ with $\mathfrak{S}(f)\subset\mathcal{O}$
and such that lifted contactomorphism $\varphi$ satisfies 
\begin{equation}
c(\psi)\leq c(\varphi).\label{eq:penultimate}
\end{equation}
This follows from Proposition \ref{prop:the_comparison_prop_prequant}: if $h_{t}$ denotes
the contact Hamiltonian of $\psi$ we choose functions $F_t:M\to\R$ supported inside $\mathcal{O}$
satisfying
\begin{equation}\label{eq:final eq}
h_t\geq F_t.
\end{equation}
The lift $\varphi$ of the corresponding path $f$ of Hamiltonian diffeomorphisms generated by $F$ satisfies the required inequality.
\end{proof}

In this setting, we can use the fact that $c_M$ satisfies the triangle inequality to obtain more information on $c$. In particular, we obtain a criterion for $c(U) $ to be finite (cf. Remark \ref{rem:when_c_is_never_zero_on_open_sets}).

\begin{Cor}
\label{cor:its_finite}
Suppose $U \subset \Sigma$ is a non-empty open set with compact closure, and suppose that $\mathrm{pr}_M(U)$ is a Hamiltonian displaceable subset of $M$. Then $c(U) <  \infty$.
\end{Cor}

\begin{proof}
We have  $c(U) \le c( \mathrm{pr}_M(U) \times S^1) = c_M( \mathrm{pr}_M(U))$, and $c_M( \mathrm{pr}_M(U)) < \infty$ by Theorem \ref{thm:HZ cap thm} below.
\end{proof} 

We also have the following result:

\begin{Prop}
\label{dis_supp_non_neg}
Suppose that $\varphi \in \widetilde{\mathrm{Cont}}_{0,c}(\Sigma ,\xi)$ has the property that $\mathrm{pr}_M(\mathfrak{S}(\varphi))$ is a Hamiltonian displaceable subset of $M$. Then $c(\varphi) \ge 0$.
\end{Prop}

\begin{proof}
First assume that $\varphi $ is the lift of an element $f \in \widetilde{\mathrm{Ham}}_c(M, d \gamma)$. The fact that $c_M(f) \ge 0$ whenever $\mathfrak{S}(f)$ is Hamiltonian displaceable is well known, but for the convenience of the reader we give the short argument here. Suppose that $g_1 \in \mathrm{Ham}_c(M ,d \gamma)$ displaces $\mathfrak{S}(f)$. Let $\{ g_t \}_{0 \le t \le 1}$ denote some path connecting $g_1$ to $\id_M$. Choose a path of paths $f^s = \{f^s_t \}_{0 \le s, t \le 1}$ connecting $f = f^1$ with the constant path $f^0_t \equiv \id_M$ such that $g_1$  displaces $\mathfrak{S}(f^s)$ for each $0 \le s \le 1$. Then we claim that 
\begin{equation}
\label{cgf=cg}
c_M(gf) = c_M(g).
\end{equation}
Indeed, the point is that any fixed point of $g_1f_1$ lies outside of $\mathfrak{S}(f)$, and hence is necessarily also a fixed point of $g_1$. The same is true if we replace $f_1$ with $f_1^s$ for any $0 \le s \le 1$, and thus it follows that $\mathrm{Spec}(\mathcal{A}_{g f^s})$ is independent of $s$. Since the function $s \mapsto c_M(g f^s)$ is continuous and $\mathrm{Spec}(\mathcal{A}_{gf^s})$ is discrete, it must be constant. This proves \eqref{cgf=cg}. 
We then argue as follows:
\begin{align}
c_M(g) & =  c_M( g f^{-1}f) \\
& \le c_M(g f^{-1}) + c_M(f) \label{tr_eq} \\
& \le c_M( g) + c_M(f) \label{cfg_app} 
\end{align}
where \eqref{tr_eq} used the triangle inequality for $c_M$ and \eqref{cfg_app} used \eqref{cgf=cg} applied to $f^{-1}$.  This implies that $c_M(f) \ge 0$. Finally to prove the general case where $\varphi$ is not necessarily the lift of a Hamiltonian path $f$, we use the same argument from the proof of Theorem \ref{thm:equality of capacities-1}. Namely, we can find a path $f$ of Hamiltonians with support inside $\mathrm{pr}_M(\mathfrak{S}(\varphi))$ such that $c_M(f) \le c(\varphi)$. Then the argument above shows that $c_M(f) \ge 0$, and hence the same is true of $c(\varphi)$.
\end{proof}

Let us quickly recall the definition of the \emph{Hofer-Zehnder capacity}.
See for instance \cite{HoferZehnder1994} for an in depth treatment.

\begin{Def}
\label{def_of_HZ}
Let $\mathcal{O}$ be an open subset of $M$. We define the \emph{Hofer-Zehnder
capacity }$c_{\textrm{HZ}}(\mathcal{O},M)$ of $\mathcal{O}$ to 
\begin{equation}
c_{\textrm{HZ}}(\mathcal{O},M):=\sup\left\{ \left\Vert H\right\Vert \mid H\mbox{ is admissible}\right\} ,
\end{equation}
where $H\in C_{c}^{\infty}(\mathcal{O},\mathbbm{R})$ is \emph{admissible
}if there exists an open set $O\subset\mathcal{O}$ such that $H|_{O}=\max H$,
and if the flow $\varphi_{H}^{t}$ has no non-constant periodic orbits
of period $\leq1$.

We also define the \emph{displacement energy }by 
\begin{equation}
e(\mathcal{O},M):=\inf\left\{ \left\Vert H\right\Vert \mid\varphi_{H}^{1}(\mathcal{O})\cap\mathcal{O}=\emptyset\right\} .
\end{equation}
The following result is due to Frauenfelder and Schlenk \cite[Corollary 8.3]{FrauenfelderSchlenk2007},
see also \cite{FrauenfelderGinzburgSchlenk2005,Schwarz2000}.\end{Def}
\begin{Thm}
\label{thm:HZ cap thm}
If $(M_1,\gamma_1)$ is a Liouville domain
then 
\begin{equation}
c_{\textrm{\emph{HZ}}}(\mathcal{O},M)\leq c_{M}(\mathcal{O})\leq e(\mathcal{O},M).
\end{equation}

\end{Thm}
Denote by $B(r)$ the open ball of radius $r$ in $\mathbbm{R}^{2m}$.
Then $c_{\textrm{HZ}}(B(r),\mathbbm{R}^{2m})=\pi r^{2}$. We can now
prove the following result, which was stated as Theorem \ref{prop:our HZ result_INTRO
} in the Introduction.
\begin{Thm}
\label{prop:our HZ result}Let $(M,d\gamma)$ denote a Liouville manifold.
Equip $\mathbbm{R}^{2m}$ with the standard symplectic form $d\lambda_{\mathrm{std}}$,
and consider the contact manifold $(\widetilde{\Sigma},\alpha+\lambda_{\mathrm{std}})$,
where $\widetilde{\Sigma}:=M\times\mathbbm{R}^{2m}\times S^{1}$. Suppose
$\mathcal{O}\subseteq M$ is open and $c_{\textrm{\emph{HZ}}}(\mathcal{O},M)<\infty$.
Choose $r_{0}>0$ such that
\begin{equation}
\left\lceil \pi r_{0}^{2}\right\rceil <\left\lceil c_{\textrm{\emph{HZ}}}(\mathcal{O},M)\right\rceil 
\end{equation}
and set 
\begin{equation}
r_{1}:=\sqrt{\tfrac{1}{\pi}c_{\textrm{\emph{HZ}}}(\mathcal{O},M)}+1
\end{equation}
Then there does not exist $\varphi\in\mathrm{Cont}_{0,c}(\widetilde{\Sigma},\alpha+\lambda_{\mathrm{std}})$
such that 
\begin{equation}
\varphi(\mathcal{O}\times B(r_{1})\times S^{1})\subset\mathcal{O}\times B(r_{0})\times S^{1}.
\end{equation}
\end{Thm}
\begin{proof}
We first prove that for $r>r_{1}$, 

\begin{equation}
c_{\textrm{HZ}}(\mathcal{O}\times B(r),M\times\mathbbm{R}^{2m})\geq c_{\textrm{HZ}}(\mathcal{O},M).\label{eq:going one way}
\end{equation}
Fix $\varepsilon>0$. We consider a cutoff function $\beta:[0,\infty)\rightarrow[0,1]$
such that $\beta(s)=1$ for $s\in[0,r-1-\varepsilon]$ and $\beta(s)=0$
for $s>r$, and such that $-1\leq\beta'(s)\leq0$ for all $s\in[0,\infty)$.
Now suppose $H$ is any admissible function on $\mathcal{O}$. Define
$H_{\beta}:M\times\mathbbm{R}^{2m}\rightarrow\mathbbm{R}$ by 
\begin{equation}
H_{\beta}(x,y):=\beta(\left|y\right|)H(x).
\end{equation}
 The symplectic gradient of $H_{\beta}$ with respect to $d\gamma\oplus d\lambda_{\mathrm{std}}$
is
\begin{equation}
X_{H_{\beta}}(x,y)=\left(\beta(\left|y\right|)X_{H}(x),H(x)X_{\beta}(y)\right).
\end{equation}
Suppose $\gamma:\mathbbm{R}\rightarrow M\times\mathbbm{R}^{2m}$ is
a non-constant periodic orbit of $X_{H_{\beta}}$, with $\gamma(t+T)=\gamma(t)$
for all $t\in\mathbbm{R}$. We shall show that $T>1$, so that $H_{\beta}$
is admissible. Write $\gamma(t)=(\gamma_{x}(t),\gamma_{y}(t))$. Then
\begin{equation}
\dot{\gamma}_{x}=\beta(\left|\gamma_{y}\right|)X_{H}(\gamma_{x}),\ \ \ \dot{\gamma}_{y}=H(\gamma_{x})X_{\beta}(\gamma_{y}).
\end{equation}
Since $\left|\beta'\right|\leq1$ we see that if $\gamma_{x}$ is
non-constant then $T>1$. But if $\gamma_{x}$ is constant, say $\gamma_{x}(t)=x_{0}$,
then we must have $H(x_{0})\ne0$. Since $\beta'$ is non-zero only
for $\left|\gamma_{y}\right|\in(r-1-\varepsilon,r)$ we necessarily
have 
\begin{equation}
T\geq\frac{1}{H(x_{0})}\pi(r-1-\varepsilon)^{2}\geq\frac{1}{c_{\textrm{HZ}}(\mathcal{O},M)}\pi(r-1-\varepsilon)^{2}.
\end{equation}
Thus as long as
\begin{equation}
\pi(r-1-\varepsilon)^{2}>c_{\textrm{HZ}}(\mathcal{O},M),\label{eq:with epsilon}
\end{equation}
$H_{\beta}$ is indeed admissible. Since clearly $\max\, H_{\beta}=\max\, H$,
we see that 
\begin{equation}
c_{\textrm{HZ}}(\mathcal{O}\times B(r),M\times\mathbbm{R}^{2n})\geq c_{\textrm{HZ}}(U,M)
\end{equation}
provided that (\ref{eq:with epsilon}) holds. Since $\varepsilon$
was arbitrary we obtain (\ref{eq:going one way}). Moreover for any
$r>0$ one always has 
\begin{equation}
e(\mathcal{O}\times B(r),M\times\mathbbm{R}^{2n})\leq\pi r^{2},
\end{equation}
as can be checked directly. The remainder of the proof is an easy
application of Theorem \ref{thm:equality of capacities-1}, Theorem
\ref{thm:HZ cap thm} and Theorem \ref{cor:non_squeezing}. Indeed,
we have 
\begin{align}
\overline{c}(\mathcal{O}\times B(r_{0})\times S^{1}) & =\left\lceil e(\mathcal{O}\times B(r_{0}),M\times\mathbbm{R}^{2m})\right\rceil \nonumber \\
 & \leq\left\lceil \pi r_{0}^{2}\right\rceil \nonumber \\
 & <\left\lceil c_{\textrm{HZ}}(\mathcal{O},M)\right\rceil \nonumber \\
 & \leq\left\lceil c_{\textrm{HZ}}(\mathcal{O}\times B(r_{1}),M\times\mathbbm{R}^{2m})\right\rceil \nonumber \\
 & \leq\overline{c}_{M\times\mathbbm{R}^{2m}}(\mathcal{O}\times B(r_{1}))\nonumber \\
 & =\overline{c}(\mathcal{O}\times B(r_{1})\times S^{1}).
\end{align}

\end{proof}
Here is an application of Theorem \ref{prop:our HZ result}, which
can be seen as a more quantitive (albeit weaker, and with more hypotheses)
version of the infinitesimal result of \cite[Theorem 1.18]{EliashbergKimPolterovich2006}.

\begin{Cor}
\label{Irie}
Suppose $X$ is a closed connected oriented Riemannian manifold which
admits a circle action $S^{1}\times X\rightarrow X$ such that the
loop $t\mapsto t\cdot p$ is not contractible for some $p\in X$.
Then if $\mathcal{O}\subset T^{*}X$ is any neighborhood of the zero
section then the conclusion of Proposition \ref{prop:our HZ result}
holds.\end{Cor}
\begin{proof}
A result of Kei Irie \cite{Irie2011} proves that in this setting
the Hofer-Zehnder capacity of the unit disc bundle $D^{*}X\subset T^{*}X$
is finite. Thus the same is true of any neighborhood $\mathcal{O}\subset T^{*}X$
of the zero section, and hence the hypotheses of Theorem \ref{prop:our HZ result}
are satisfied.
\end{proof}

\appendix

\section{The ``compactly supported Reeb
flow''}
\label{app:A-'compactly-supported}

In this Appendix we continue to work in the setting from the previous section.  Thus $\Sigma = M \times S^1$ is a prequantisation space associated to the completion of a Liouville domain $(M_1, d \gamma_1)$. Our aim is to construct a ``compactly supported Reeb flow'' whose
support is contained in a tubular neighborhood of a closed Reeb orbit, and explicitly compute the spectral value. This result has been alluded to in Remarks \ref{rem:reversing_prev_cor} and \ref{rem:cptly_sprtd_rf}.

\begin{Thm}
\label{prop:compactly supported reeb}Suppose $(\Sigma = M \times S^1,\xi)$ satisfies
Assumption \textbf{\emph{(C)}}. Let $x(t) = (y_0,t)$ denote a closed embedded Reeb orbit (for some fixed $y_0 \in M$.) Then there exists $\rho_{0}>0$  and a neighborhood $B$ of $y_0 $ in $M$ with the following significance: For
all $\rho\in\mathbbm{R}$ with $| \rho |<\rho_{0}$, there
exists an exact contactomorphism $\vartheta^{\rho}\in\widetilde{\mathrm{Cont}}_{0}(\Sigma,\xi)$
with $\mathfrak{S}(\vartheta^{\rho})\subset B \times S^1$ with the property
that if $x\in\mathfrak{S}(\vartheta^{\rho})$ is a translated point
with of $\vartheta^{\rho}$ then 
\begin{equation}
\vartheta^{\rho}(x)=\theta^{\rho}(x).
\end{equation}
In other words, from the point of view of translated points, $\vartheta^{\rho}$
is ``the Reeb flow supported on $x$''. Moreover if $B' \subset B $
is any neigborhood of $y_0$ then for $| \rho |$
sufficiently small we have $\mathfrak{S}(\vartheta^{\rho})\subset B' \times S^1$. 
The spectral value $c(\vartheta^{\rho})$ is given by 
\begin{equation}
c(\vartheta^{\rho})=\begin{cases}
0, & 0\leq\rho<\rho_{0},\\
-\rho, & -\rho_{0}<\rho\leq0.
\end{cases}
\end{equation}
\end{Thm}

\textbf{\underline{Convention:}} In this appendix we equip $\mathbbm{R}^{2n}\backslash\{0\}$
with polar coordinates $(s,\phi)$ where $s\in(0,\infty)$ and $\phi=(\phi_{1},\dots,\phi_{2n-1})$
with $\phi_{j}\in\mathbbm{R}/2\pi\mathbbm{Z}$. In these coordinates the standard contact form $\alpha_{\mathrm{std}}$
is given by
\begin{equation}
\label{eq:alpha_std}
\alpha_{\mathrm{std}}=\sum_{j}\tfrac{1}{2}s^{2}d\phi_{j}+d\tau.
\end{equation}
This has the slightly unfortunate consequence that $\tau$ is 1-periodic but the $\phi_{j}$ are $2\pi$-periodic! These conventions are chosen so that $c_{\mathbbm{R}^{2n}}(B(r))=\pi r^{2}$ instead of  $\tfrac{1}{2}r^{2}$. 

 \begin{proof}[Proof of Theorem \ref{prop:compactly supported reeb}]
The argument is local in $M$, and hence it is sufficient to prove the result in the special case $M = \R^{2n}$.  Thus $\Sigma = \R^{2n} \times S^1$ and $\alpha = \alpha_{\mathrm{std}}$ is given by \eqref{eq:alpha_std}. The Reeb vector field $R$ of $\alpha$ is just $\frac{\partial}{\partial\tau}$,
and the Reeb flow $\theta^{t}$ is given by 
\begin{equation}
\theta^{t}(s,\phi,\tau)=(s,\phi,\underset{\textrm{mod 1}}{\underbrace{\tau+t}}).
\end{equation}
Fix $\rho\in\mathbbm{R}$ such that $0< | \rho |<\pi r^{2}$.
Let $f:[0,\infty)\times[0,\infty)\rightarrow\mathbbm{R}$ denote a
smooth function with the following properties: 
\begin{enumerate}
\item There exists $\varepsilon>0$ such that $f(s)=\rho$ for $0\leq s\leq\varepsilon$
and $f(s)=0$ for $r-\varepsilon\leq s\leq r$.
\item If $\rho<0$ then $f'(s)\geq0$ for all $s$. If $\rho>0$ then $f'(s)\leq0$
for all $s$.
\item If $\rho<0$ then $2\pi s-f'(s)>0$ for all $s>0$. If $\rho<0$ then
$2\pi s+f'(s)<0$ for all $s>0$.
\end{enumerate}
Note that such a function only exists because $\left|\rho\right|<\pi r^{2}$.
Indeed, if $\rho<0$ then since $2\pi s-f'(s)>0$ one has
\begin{equation}
-\rho=\int_{0}^{r}f'(s)ds<\int_{0}^{r}2\pi sds=\pi r^{2}.\label{eq:its empty}
\end{equation}
Conversely it is easy to see that when $\left|\rho\right|<\pi r^{2}$
such functions really do exist. Now consider the contactomorphism
$\vartheta^{\rho}$ of $\mathbbm{R}^{2n}\times S^{1}$ whose contact
Hamiltonian $h_{t}:\mathbbm{R}^{2n}\times S^{1}$ is given by 
\begin{equation}
h_{t}(s,\zeta,\tau)=f(r).
\end{equation}
The contact vector field $X_{t}$ of $h_{t}$ is defined by the equations
\begin{equation}
\alpha(X_{t})=h_{t},\ \ \ i_{X_{t}}d\alpha=dh_{t}(R)\alpha-dh_{t}.
\end{equation}
This gives
\begin{equation}
X_{t}(s,\phi,\tau)=\sum_{j}\frac{f'(s)}{s}\frac{\partial}{\partial\phi_{j}}+\left(f(s)-\frac{sf'(s)}{2}\right)\frac{\partial}{\partial\tau}.
\end{equation}
We can integrate this to obtain 
\begin{equation}
\vartheta_{t}^{\rho}(s,\phi,\tau)=\left(s,\phi_{1}+\frac{f'(s)}{s}t,\dots,\phi_{2n-1}+\frac{f'(s)}{2}t,\tau+\left(f(s)-\frac{sf'(s)}{2}\right)t\right),
\end{equation}
and hence translated points of $\vartheta_{1}^{\rho}$ are tuples
$(s,\phi,\tau)$ with 
\begin{equation}
\frac{f'(s)}{s}\in2\pi\mathbbm{Z},\label{eq:tp1}
\end{equation}
and the time-shift is given by 
\begin{equation}
\eta=f(s)-\frac{sf'(s)}{2}.\label{eq:tp2}
\end{equation}
By assumption one never has $f'(s)/2\pi s\in\mathbbm{Z}$ unless $f'(s)=0$.
In other words, translated points only occur when $0\leq s\leq\varepsilon$
or when $ $$r-\varepsilon\leq s\leq\infty$. In particular, the only
translated points of $\vartheta^{\rho}$ that lie in the interior
of the support of $\vartheta^{\rho}$ are the points in $B(\varepsilon)\times S^{1}$.
Since $\vartheta^{\rho}=\theta^{\rho}$ on $B(\varepsilon)\times S^{1}$,
this justifies our claim that `from the point of view of translated
points', $\vartheta^{\rho}$ is the Reeb flow. 

To complete the proof let us compute the spectral value of $\vartheta^{\rho}$.
Note that the \emph{contractible }action spectrum of $\mathcal{A}_{\vartheta^{\rho}}$
is just $\{0,-\rho\}$, and hence we certainly have $c(\vartheta^{\rho})\in\{0,-\rho\}$.
For $\rho<0$, one has $h_{t}<0$ on the interior of its support and
hence by Corollary \ref{cor:posi_but_not_negi_prequant} one has $c(\vartheta^{\rho})>0$,
which implies $c(\vartheta^{\rho})=-\rho$. Thus for $\rho>0$ we must have $c(\vartheta^{\rho})=0$.
This completes the proof.
\end{proof}

\bibliographystyle{amsalpha}
\bibliography{willmacbibtex}
\end{document}